\documentclass[12pt]{article}

\usepackage{amsmath,amssymb,amsfonts,amsthm,hyperref,bbm,latexsym,epsfig}
\usepackage{graphicx,multirow} 
\usepackage{caption}

\newcommand{\bs}{\begin{subequations}}
\newcommand{\es}{\end{subequations}}
\newcommand{\be}{\begin{equation}}
\newcommand{\ee}{\end{equation}}
\newcommand{\ba}{\begin{eqnarray}}
\newcommand{\ea}{\end{eqnarray}}

\allowdisplaybreaks

\begin{document}

\title{
\normalsize \hfill CFTP/17-001
\\[6mm]
\LARGE {\tt GAP} listing of the finite subgroups of $U(3)$ \\
of order smaller than 2\,000}

\author{
\addtocounter{footnote}{2}
Darius~Jur\v{c}iukonis$^{(1)}$\thanks{E-mail: \tt darius.jurciukonis@tfai.vu.lt}
\ {\normalsize and}
Lu\'\i s~Lavoura$^{(2)}$\thanks{E-mail: \tt balio@cftp.tecnico.ulisboa.pt}
\\*[3mm]
$^{(1)} \! $
\small Vilnius University, 
\small Institute of Theoretical Physics and Astronomy, \\
\small Saul\.etekio ave.\ 3, LT-10222 Vilnius, Lithuania 
\\[2mm]
$^{(2)} \! $
\small Universidade de Lisboa, Instituto Superior T\'ecnico, CFTP, \\
\small 1049-001 Lisboa, Portugal
\\*[2mm]
}

\date{\today}

\maketitle

\begin{abstract}
  We have sorted the {\tt SmallGroups} library
  of all the finite groups of order smaller than 2\,000
  to identify the groups that
  possess a faithful three-dimensional irreducible representation
  (`irrep')
  and cannot be written as the direct product
  of a smaller group times a cyclic group.
  Using the computer algebra system {\tt GAP},
  we have scanned all the three-dimensional irreps of each of those groups
  to identify those that are subgroups of $SU(3)$;
  we have labelled each of those subgroups of $SU(3)$
  by using the extant complete classification
  of the finite subgroups of $SU(3)$.
  Turning to the subgroups of $U(3)$ that are not subgroups of $SU(3)$,
  we have found the generators of all of them
  and classified most of them in series according to their generators
  and structure.
\end{abstract}

\newpage

\section{Introduction}

Many high-energy physicists are thrilled by the prospect that
the numerical entries of the leptonic mixing matrix (PMNS matrix)
might be related to some small (or maybe not so small) finite group.
Many specific finite groups have been considered,
like for instance $A_4$~\cite{a4},
$S_4$~\cite{s4},
$S_3$~\cite{s3},
$T_7$~\cite{t7},
$A_5$~\cite{a5},
$\Delta(27)$~\cite{delta27},
the group series $\Delta \left( 6 n^2 \right)$~\cite{delta6n2},
the groups $\Sigma \left( n \varphi \right)$~\cite{Sigmanphi},
and so on.
Most of the finite groups considered are subgroups of $SU(3)$;
those subgroups are especially inviting
because a complete classification of them,
and their generators,
have been known for over a century~\cite{su3classification}.
On the contrary,
there is no complete classification of the
finite subgroups of $U(3)$,\footnote{In this paper,
  whenever we use the expression ``finite subgroups of $U(3)$''
  we usually mean only the subgroups of $U(3)$
  that are not subgroups of $SU(3)$.}
though a few series of those subgroups have been derived in ref.~\cite{ludl}.
At least one finite subgroup of $U(3)$
has already been utilized in particle physics~\cite{sigma81}.

Although a full theoretical study of each individual group
can always be undertaken,
for large groups such a study becomes impractical
and it is convenient to have recourse to the computer algebra system {\tt GAP},
which is tailored to deal with finite groups
and can readily furnish the structure,
irreducible representations (`irreps'),
character table,
and so on,
of each of them.
{\tt GAP} is supplemented by the {\tt SmallGroups} library,
which contains,
in particular,
all the finite groups of order smaller than 2\,000.
In that library each finite group has an identifier
$\left[ o, j \right]$,
where $o \ge 1$ is the order,
\textit{i.e.}\ the number of elements,
of the group and $j \ge 1$ is an integer which distinguishes
among the non-isomorphic groups of identical order.
For instance,
the group with {\tt SmallGroups} identifier $\left[ 4, 1 \right]$
is the cyclic group\footnote{{\tt SmallGroups} uses $C_n$
  to denote the cyclic group of order $n$,
  instead of the more usual notation $\mathbbm{Z}_n$.
{\tt SmallGroups} uses the notation $E(n)$ for the $n$'th root of unity.}
$\mathbbm{Z}_4 \cong \left\{ 1,\, i,\, -1,\, -i \right\}$
while the group with {\tt SmallGroups} identifier $\left[ 4, 2 \right]$
is the direct product of cyclic groups
$\mathbbm{Z}_2 \times \mathbbm{Z}_2 \cong \left\{
\left( 1, 1 \right),\, \left( 1, -1 \right),\, \left( -1, 1 \right),\,
\left( -1, -1 \right) \right\}$;
{\tt SmallGroups} informs us that there are,
in fact,
only these two non-isomorphic groups with four elements.
A
{\tt SmallGroups} listing of all the finite groups
of order up to 100,
together with their structure,\footnote{{\tt SmallGroups} informs us about the
  structure of each group.
  This is given in terms of direct products (denoted `$\times$'),
  semi-direct products (denoted `$\rtimes$'),
  or group extensions (denoted `$.$').
  A pedagogical explanation of these concepts may be found,
  for instance,
  in ref.~\cite{vdovin}.}
was published in ref.~\cite{parattu}.
A
{\tt SmallGroups} listing of the finite groups of order up to 512
that have a faithful three-dimensional irrep
and are not the direct product of a cyclic group and some other group
was published in ref.~\cite{ludl}.

However,
{\tt SmallGroups}
lists
the groups of the same order
in a way that does not allow one to extract much information on them.
For instance,
\begin{description}
\item the group $\left[ 12, 3 \right] \cong A_4$ is a subgroup of $SU(3)$
  and has a three-dimensional faithful irrep;
\item the groups $\left[ 12, 1 \right]$ and $ \left[ 12, 4 \right] \cong D_6$
  are subgroups of $SU(3)$ but do not possess three-dimensional irreps;
\item the group $\left[ 12, 2 \right] \cong \mathbbm{Z}_{12}$
  is a subgroup of
  $U(1) \subset U(3)$;
\item the group $\left[ 12, 5 \right] \cong
  \mathbbm{Z}_6 \times \mathbbm{Z}_2$
  is a subgroup of $U(1) \times U(1)$ but not of $U(3)$.
\end{description}
%
The first step in this work was to survey the whole {\tt SmallGroups} list
of groups of order smaller than 2\,000 in order to identify the ones
that
\begin{description}
  \item have at least one faithful three-dimensional irreducible
    representation;
  \item cannot be written as the direct product of a smaller group
    and a cyclic group.
\end{description}
The second step in this work was to pick each of the groups above
and ask {\tt GAP} to compute the determinant of each of the matrices
in each of its three-dimensional representations.
If there is a three-dimensional representation in which all the matrices
have unit determinant,
then the group is a subgroup of $SU(3)$;
otherwise the group is not a subgroup of $SU(3)$
but it is a subgroup of $U(3)$---because every representation of a finite group
is equivalent to a representation through unitary matrices.
In this way,
we have separated the subgroups of $SU(3)$ from the subgroups of $U(3)$.

A complete classification of all the finite subgroups of $SU(3)$
has long existed~\cite{su3classification}.
There are groups (so-called type A) of diagonal matrices,
\textit{i.e.}\ Abelian groups;
they may be written as direct products
of cyclic factors and do not concern us here.
Then there are the subgroups of $U(2)$,
which are called type B;
their three-dimensional representations are
(just as the ones of type~A subgroups)
reducible and therefore they also do not concern us.
Of interest to us are the type C and type D groups,
which were best characterized in ref.~\cite{grimus},
and also the `exceptional' groups.
In this work we give the {\tt SmallGroups} identifiers
of all the $SU(3)$ subgroups of types C and D,
together with their classification according to ref.~\cite{grimus},
and also the {\tt SmallGroups} identifiers of the exceptional subgroups.
This is done in section~\ref{su3}.

There is no theoretical classification
of all the finite subgroups of $U(3)$.
We feel that having a complete listing of all those subgroups
of order less than 2\,000,
together with their generators,
may be a useful step towards achieving such a classification;
at the very least,
it allows one to get a feeling for what it might look like.
Therefore,
in this work we give the {\tt SmallGroups} identifiers
of all the finite $U(3)$ subgroups,
together with their generators.
We also partially unite those subgroups in series,
\textit{viz.}\ in sets of groups that have related generators depending on one,
two,
or sometimes three integers.
This is done in section~\ref{u3}.

We also give,
for every finite subgroup of $U(3)$,
the dimensions of all its inequivalent irreps,
as determined by {\tt GAP}.

In section~\ref{proc} we explain our procedure.
In an appendix we provide tables
of all the finite subgroups of $U(3)$
that have a faithful three-dimensional irrep
and are not isomorphic to the direct product of a smaller group
and a cyclic group.
We give separate tables
for the groups that are subgroups of $SU(3)$
and for the groups that are not subgroups of $SU(3)$.
In those tables,
we order the groups according to their {\tt SmallGroups} classification,
\textit{viz.}\ in increasing order first of $o$ and then of $j$
in their $\left[ o, j \right]$ identifiers.

\section{{\tt GAP} procedures} \label{proc}

{\tt GAP}~\cite{GAP}
is a computer algebra system that provides a programming language,
including many functions that implement algebraic algorithms.
It is supplemented by many libraries containing a large amount of data
on algebraic objects.
Using {\tt GAP} it is possible to study groups and their representations,
display the character tables,
find the subgroups of larger groups,
identify groups given through their generating matrices,
and so on.

{\tt GAP} allows access to the {\tt SmallGroups} library
through the {\tt SmallGroups} package~\cite{SmallGroups}.
That library contains \emph{all}\/ the finite groups of `small'
orders,\footnote{The order of a finite group is the number of its elements.}
\textit{viz.}\ less than a certain upper bound
and also orders whose prime factorization is small in some sense.
The groups are ordered by their orders;
for each of the available orders,
a complete list of non-isomorphic groups is given.
{\tt SmallGroups} contains all the groups of order less than 2\,000
except order 1024,
because there are many thousands of millions of groups of order 1024.
{\tt SmallGroups} also contains other groups with some specific orders
larger than 2\,000.

The {\tt SmallGroups} library has an identification function
which returns the {\tt SmallGroups} identifier of any given group.
For each generic group in the library
there are effective recognition algorithms available.
To identify encoded and insoluble groups,
two approaches are used:
one is a general algorithm to solve the isomorphism problem
for $p$-groups,\footnote{A $p$-group,
  where $p$ is a prime number,
  is a group in which each element has a power of $p$ as its order.
  That is,
  for each element $g$ of a $p$-group,
  there is a non-negative integer $n$ such that
  the product of $p^n$ copies of $g$,
  and not less,
  is equal to the identity element $e$.
  (But,
  the integer $n$ is in general different for different elements $g$
  of the group.)}
the second one uses the invariants\footnote{In the {\tt SmallGroups} library
there is a list of distinguishing invariants
for all encoded groups except those of orders 512 and 1536.
This list of invariants is compressed.
It provides an efficient approach to identify
any encoded group in the library.}
of stored groups~\cite{besche}.
Using these methods,
it is possible to identify all the groups in the library,
except for orders 512,
1536,
and some orders above 2\,000.
For the identification of groups we use {\tt GAP} command
\be {\tt IdGroup(.)}. \ee

In our work,
firstly we have scanned the {\tt SmallGroups} library
and extracted therefrom all the groups with three-dimensional irreps.
Using the {\tt GAP} command
\be {\tt G:=SmallGroup([o,j])}, \ee
one lets $G$ denote the group with identifer $\left[ o, j \right]$
in the {\tt SmallGroups} library.
The command
\be {\tt NumberSmallGroups(o)} \ee
allows one to find out how many groups there are for a chosen order $o$
and thus automates the scanning of library.
For a given group $G$,
{\tt GAP} offers the possibility
to calculate the irreps by using the command
\be {\tt repG:=IrreducibleRepresentations(G)}. \ee
It is possible to display all the irreps by using the {\tt GAP} command
\be {\tt Display(CharacterTable(G))} \ee
too;
however,
the labeling of the irreps may differ from the labeling
received through the command
\be {\tt IrreducibleRepresentations(G)}. \ee
It is convenient to select all the three-dimensional irreps
by using the command
\be {\tt
  repG3:=Filtered(repG,x\mbox{--}\!>\!Length(Identity(G)\verb!^!x)=3)}. \ee
One may select all the elements 
of a given group $G$ through the command
\be {\tt elG:=Elements(G)}. \ee
Then,
the command
\be \label{biughp}
    {\tt elGlist:=List(elG,x\mbox{--}\!>\!x\verb!^!repG3[i])}, \ee
%
where the integer $i$ parameterizes the loop,
allows one to list all the elements of the chosen irrep.
We have selected the groups from the {\tt SmallGroups} library
that have at least one faithful\footnote{In order to identify
  the faithful irreps,
  we have compared all the matrices in each irrep.
  If different elements of the group are represented by different matrices
  in the irrep,
  then the irrep is faithful.}
three-dimensional irrep.
Then,
by using the {\tt GAP} command that gives the structure of a group,
\textit{viz.}
\be \label{bvury} {\tt StructureDescription(G)}, \ee
we have discarded
the groups that are direct products with a cyclic group. 

There are 10\,494\,213 groups of order 512
and 408\,641\,062 groups of order 1536.
However,
the groups of order 512 do not possess three-dimensional irreps
because 512 cannot be divided by three,
therefore we did not need to consider them.
On the other hand,
the number of groups of order 1536 is too large
for all of them to be scanned in the way described above.
Therefore,
we have used the conjecture in ref.~\cite{chinese} that
both nilpotent groups
and groups with a normal Sylow 3-subgroup\footnote{These two concepts
of group theory have been explained in ref.~\cite{rows}.}
do not have three-dimensional faithful irreps.
Utilizing the command
\be {\tt SmallGroupsInformation(o)} \ee
one gets the information about the arrangement of the groups of a given order.
Using this information,
we have determined the scanning range of groups of order 1536.
To check whether the group is nilpotent,
the command
\be {\tt IsNilpotentGroup(G)} \ee
may be used,
while
\be {\tt NilpotencyClassOfGroup(G)} \ee
gives the nilpotency class of the group $G$.
The Sylow 3-subgroups of a group $G$
may be found by typing the command
\be {\tt SylowSubgroup(G,3)}. \ee
We have found that only four groups of order 1536
have faithful three-dimensional irreps
and cannot be written as the direct product of a smaller group
and a cyclic group.

For groups that have faithful three-dimensional irreps,
we have asked {\tt GAP} to compute the determinant
of each of the matrices in each of its three-dimensional representations.
This was done through the command
\be {\tt DeterminantMat(elGlist[i])}. \ee
If there is a three-dimensional representation
in which all the matrices have unit determinant,
then the group is a subgroup of $SU(3)$;
if there is no such representation,
then the group is not a subgroup of $SU(3)$,
but it is a subgroup of $U(3)$
because it has a three-dimensional representation
and because all the representations of finite groups are equivalent
to representations through unitary matrices.

We have used different methods
in order to classify the groups in the lists of the subgroups
of $U(3)$ and $SU(3)$.
One of the methods is the analysis of the generators
of the three-dimensional irreps.
The command
\be {\tt genG:=GeneratorsOfGroup(G)} \ee
returns a list of generators of the group $G$.
The generators of the three-dimensional irreps
may be listed through the command
\be {\tt List(genG,x\mbox{-}\!>\!x\verb!^!repG3[i])}. \ee
By looking at these lists we have tried to find regularities in the generators.
Another strategy was looking at the structures of the groups
and sorting groups with analogous structures.

When one has some generators,
say three matrices $M1$,
$M2$,
and $M3$,
a group $G$ may be generated through the command
\be {\tt G:=Group([M1, M2, M3])}. \ee
Afterwards this group may be identified by finding its order,
using the command
\be {\tt Order(G)} \ee
or by counting the elements of the group through
\be {\tt Size(elG)}. \ee
Afterwards one may discover the {\tt SmallGroups} identifier of $G$
by using the command
\be {\tt IdGroup(G)}. \ee
The identification of some groups with large order
may require a long computational time,
therefore some hints about the group classification
may be acquired by analyzing the group structure---using
the command~\eqref{bvury}---or by comparing
the traces of the group matrices,
determined through the command
\be {\tt List(elG,x\mbox{--}\!>\!Trace(x))}. \ee

\section{Finite subgroups of $SU(3)$} \label{su3}

In this section we give the generators and the {\tt SmallGroups} identifiers
of all the finite subgroups of $SU(3)$ that
\begin{itemize}
\item have a faithful three-dimensional irrep,
\item cannot be written as the direct product of a smaller group
  and a cyclic group,
\item have less than 2\,000 elements.
\end{itemize}

\subsection{Generators}
  
We firstly define a few $3 \times 3$ matrices that act as generators
of the various $SU(3)$ subgroups.
All those matrices have,
of course,
unit determinant.

The matrices
\bs
\ba
\label{E}
E &\equiv& \left( \begin{array}{ccc} 0 & 1 & 0 \\ 0 & 0 & 1 \\ 1 & 0 & 0
\end{array} \right),
\\
\label{Igen}
I &\equiv&
\left( \begin{array}{ccc} 0 & 0 & -1 \\ 0 & -1 & 0 \\ -1 & 0 & 0
  \end{array} \right)
\ea
\es
are especially useful.
Let $n \ge 1$ be an integer.
Then,
\be
\label{nu}
L_n \equiv \mathrm{diag} \left( 1,\, \nu,\, \nu^{-1} \right),
\quad \mbox{where}\ \nu = \exp{\left( 2 i \pi / n \right)}.
\ee
Let $n \ge 1$ and $k \ge 1$ be integers.
We define
\be
\label{Bmk}
B_{n,k} \equiv \mathrm{diag} \left( \nu,\, \nu^k,\, \nu^{-1-k} \right),
\quad \mbox{where}\ \nu = \exp{\left( 2 i \pi / n \right)}.
\ee
Let $n \ge 1$ and $r \ge 1$ be integers.
We define
\be
\label{G}
G_{n,r} \equiv \mathrm{diag} \left( 1,\, \nu^{-r},\, \nu^r \right),
\quad \mbox{where}\ \nu = \exp{\left( 2 i \pi / n \right)},
\ee
%
\textit{i.e.}\ $G_{n,r} = \left( L_n \right)^{-r}$.

\subsection{The groups $\Delta \left( 3 n^2 \right)$
  and $\Delta \left( 6 n^2 \right)$}

For $n \ge 1$,
the groups $\Delta \left( 3 n^2 \right)$ have structure
$\left( \mathbbm{Z}_n \times \mathbbm{Z}_n \right) \rtimes \mathbbm{Z}_3$
and order $3 n^2$;\footnote{We adopt the convention that $\mathbbm{Z}_1$
  is the trivial group,
  \textit{i.e.}\ the group that has only one element,
  \textit{viz.}\ the identity element $e$.}
the groups $\Delta \left( 6 n^2 \right)$ have structure
$\left[ \left( \mathbbm{Z}_n \times \mathbbm{Z}_n \right) \rtimes \mathbbm{Z}_3
  \right] \rtimes \mathbbm{Z}_2$ and order $6 n^2$.
The group $\Delta \left( 3 n^2 \right)$ is generated by the matrices $E$
and $L_n$;
the group $\Delta \left( 6 n^2 \right)$ is generated by the matrices $E$,
$I$,
and $L_n$.
The {\tt SmallGroups} identifiers
of the groups $\Delta \left( 3 n^2 \right)$ of order smaller than 2\,000
are given in table~\ref{delta3n2};\footnote{The group
  $\Delta \left( 3\times1^2 \right) \cong \mathbbm{Z}_3
  \cong \left[ 3, 1 \right]$
  is not included in table~\ref{delta3n2} because it is a cyclic group.}
\begin{table}
\begin{center}
\renewcommand{\arraystretch}{1.4}
\begin{tabular}{|c||c|c|c|c|c|c|} \hline
  $n$ & 2 & 3 & 4 & 5 \\*[-2mm]
  identifier &
  $\left[ 12, 3 \right]$ &
  $\left[ 27, 3 \right]$ &
  $\left[ 48, 3 \right]$ &
  $\left[ 75, 2 \right]$ \\ \hline
  $n$ & 6 & 7 & 8 & 9 \\*[-2mm]
  identifier &
  $\left[ 108, 22 \right]$ &
  $\left[ 147, 5 \right]$ &
  $\left[ 192, 3 \right]$ &
  $\left[ 243, 26 \right]$ \\ \hline
  $n$ & 10 & 11 & 12 & 13 \\*[-2mm]
  identifier &
  $\left[ 300, 43 \right]$ &
  $\left[ 363, 2 \right]$ &
  $\left[ 432, 103 \right]$ &
  $\left[ 507, 5 \right]$ \\ \hline
  $n$ & 14 & 15 & 16 & 17 \\*[-2mm]
  identifier &
  $\left[ 588, 60 \right]$ &
  $\left[ 675, 12 \right]$ &
  $\left[ 768, 1083477 \right]$ &
  $\left[ 867, 2 \right]$ \\ \hline
  $n$ & 18 & 19 & 20 & 21 \\*[-2mm]
  identifier &
  $\left[ 972, 122 \right]$ &
  $\left[ 1083, 5 \right]$ &
  $\left[ 1200, 384 \right]$ &
  $\left[ 1323, 43 \right]$ \\ \hline
  $n$ & 22 & 23 & 24 & 25 \\*[-2mm]
  identifier &
  $\left[ 1452, 34 \right]$ &
  $\left[ 1587, 2 \right]$ &
  $\left[ 1728, 1291 \right]$ &
  $\left[ 1875, 16 \right]$ \\ \hline
\end{tabular}
\end{center}
\captionsetup{width=11.0cm}
\caption{The {\tt SmallGroups} identifiers
  of the groups $\Delta \left( 3 n^2 \right)$ with order smaller than 2\,000.
\label{delta3n2}}
\end{table}
the {\tt SmallGroups} identifiers
of the groups $\Delta \left( 6 n^2 \right)$ of order smaller than 2\,000
are given in table~\ref{delta6n2}.\footnote{The group
  $\Delta \left( 6\times1^2 \right) \cong S_3
  \cong \left[ 6, 1 \right]$ is not included in table~\ref{delta6n2}
  because its three-dimensional representations are reducible.}
\begin{table}
\begin{center}
\renewcommand{\arraystretch}{1.4}
\begin{tabular}{|c||c|c|c|c|c|c|} \hline
  $n$ & 2 & 3 & 4 & 5 \\*[-2mm]
  identifier &
  $\left[ 24, 12 \right]$ &
  $\left[ 54, 8 \right]$ &
  $\left[ 96, 64 \right]$ &
  $\left[ 150, 5 \right]$ \\ \hline
  $n$ & 6 & 7 & 8 & 9 \\*[-2mm]
  identifier &
  $\left[ 216, 95 \right]$ &
  $\left[ 294, 7 \right]$ &
  $\left[ 384, 568 \right]$ &
  $\left[ 486, 61 \right]$ \\ \hline
  $n$ & 10 & 11 & 12 & 13 \\*[-2mm]
  identifier &
  $\left[ 600, 179 \right]$ &
  $\left[ 726, 5 \right]$ &
  $\left[ 864, 701 \right]$ &
  $\left[ 1014, 7 \right]$ \\ \hline
  $n$ & 14 & 15 & 16 & 17 \\*[-2mm]
  identifier &
  $\left[ 1176, 243 \right]$ &
  $\left[ 1350, 46 \right]$ &
  $\left[ 1536, 408544632 \right]$ &
  $\left[ 1734, 5 \right]$ \\ \hline
  $n$ & 18 & & & \\*[-2mm]
  identifier &
  $\left[ 1944, 849 \right]$ & & & \\ \hline
\end{tabular}
\end{center}
\captionsetup{width=11.7cm}
\caption{The {\tt SmallGroups} identifiers
  of the groups $\Delta \left( 6 n^2 \right)$ with order smaller than 2\,000.
\label{delta6n2}}
\end{table}

The group $\Delta \left( 3\times2^2 \right)$ is isomorphic to $A_4$,
the group of the even permutations of four objects,
and also to the symmetry group of the regular tetrahedron.
The group $\Delta \left( 6\times2^2 \right)$ is isomorphic to $S_4$,
the group of all the permutations of four objects,
and also to the symmetry group of the cube and of the regular octahedron.

When $n$ cannot be divided by three,
the group $\Delta \left( 3 n^2 \right)$ has three singlet irreps
and $\left. \left( n^2 - 1 \right) \right/ 3$ triplet irreps;
when $n$ is a multiple of three,
$\Delta \left( 3 n^2 \right)$ has nine inequivalent singlet irreps
and $n^2/3 - 1$ inequivalent triplet irreps.
The group $\Delta \left( 6 n^2 \right)$
has~\cite{delta6n2,su3classification,finite},
for any $n \ge 2$,
two inequivalent singlet irreps and $2(n-1)$ inequivalent triplet irreps.
When $n$ is not a multiple of three,
$\Delta \left( 6 n^2 \right)$ has one doublet irrep
and $(n-1)(n-2)/6$ six-dimensional irreps;
when $n$ is a multiple of three,
$\Delta \left( 6 n^2 \right)$ has four inequivalent doublet irreps
and $n(n-3)/6$ six-dimensional irreps.

\subsection{The groups $C_{n,l}^{(k)}$}

We use the notation of ref.~\cite{grimus}.
The groups $C_{n,l}^{(k)}$ have structure
$\left( \mathbbm{Z}_n \times \mathbbm{Z}_l \right) \rtimes \mathbbm{Z}_3$
and order $3 n l$.
The integer $l$ is positive.
The integer $n$ may be written $n = r l$,
where $r$ is another positive integer.
The integer $r$ may be either
\begin{enumerate}
  \item a product of prime numbers $p_1, p_2, \ldots$
which are of the form $p_j = 6 i_j + 1$,
where the numbers $i_j$ are integers, or \label{1}
\item three times a product of prime numbers as in~\ref{1}. \label{2}
  \end{enumerate}
In case~\ref{1},
$l$ may be any positive integer;
in case~\ref{2},
$l$ must be a multiple of three.
The integer $k$ is a function of $r$ defined by $1 + k + k^2 = 0\ \mbox{mod}\ r$
and $k \le \left. \left( r - 1 \right) \right/ 2$.
For most values of $r$ there is only one possible $k$,
but for some $r$ more than one (usually two) $k$ are possible.
The values of $r$,
$k$,
and $l$ that produce groups $C_{n,l}^{(k)}$ with order smaller than 2\,000
are given in tables~\ref{rmn} and~\ref{rmn2}.
\begin{table}
\begin{center}
\renewcommand{\arraystretch}{1.2}
\begin{tabular}{|c||c|c|c|c|c|c|c|} \hline
  $r$ & 3 & 7 & 13 & 19 & 21 & 31 & 37 \\*[-1mm]
  $k$ & 1 & 2 & 3 & 7 & 4 & 5 & 10 \\*[-1mm]
  $l$ & 3, 6, 9, 12 & 1 to 9 & 1 to 7 & 1 to 5 & 3 & 1 to 4 & 1 to 4 \\ \hline
  $r$ & 39 & 43 & 49 & 57 & 61 & 67 & 73 \\*[-1mm]
  $k$ & 16 & 6 & 18 & 7 & 13 & 29 & 8 \\*[-1mm]
  $l$ & 3 & 1, 2, 3 & 1, 2, 3 & 3 & 1, 2, 3 & 1, 2, 3 & 1, 2, 3 \\ \hline
  $r$ & 79 & 91 & 97 & 103 & 109 & 127 & 133 \\*[-1mm]
  $k$ & 23 & 9, 16 & 35 & 46 & 45 & 19 & 11,\ 30 \\*[-1mm]
  $l$ & 1, 2 & 1, 2 & 1, 2 & 1, 2 & 1, 2 & 1, 2 & 1, 2 \\ \hline
  $r$ & 139 & 151 & 157 & 163 & 169 & 181 & \\*[-1mm]
  $k$ & 42 & 32 & 12 & 58 & 22 & 48 & \\*[-1mm]
  $l$ & 1, 2 & 1, 2 & 1, 2 & 1, 2 & 1 & 1 & \\ \hline
\end{tabular}
\end{center}
\captionsetup{width=12.1cm}
\caption{The lowest values of $r$,
  and the corresponding values of $k$ and $l$,
  that produce groups $C_{rl,l}^{(k)}$ with order $3 r l^2 < 2\,000$.
\label{rmn}}
\end{table}
\begin{table}
\begin{center}
\renewcommand{\arraystretch}{1.2}
\begin{tabular}{|c||c|c|c|c|c|c|} \hline
$r$ & 193 & 199 & 211 & 217 & 223 & 229 \\*[-1mm]
$k$ & 84 & 92 & 14 & 25,\ 67 & 39 & 94 \\ \hline
  $r$ & 241 & 247 & 259 & 271 & 277 & 283 \\*[-1mm]
  $k$ & 15 & 68, 87 & 100, 121 & 28 & 116 & 44 \\ \hline 
  $r$ & 301 & 307 & 313 & 331 & 337 & 343 \\*[-1mm]
  $k$ & 79, 135 & 17 & 98 & 31 & 128 & 18 \\ \hline 
  $r$ & 349 & 361 & 367 & 373 & 379 & 397 \\*[-1mm]
  $k$ & 122 & 68 & 83 & 88 & 51 & 34 \\ \hline
  $r$ & 403 & 409 & 421 & 427 & 433 & 439 \\*[-1mm]
  $k$ & 87, 191 & 53 & 20 & 74, 135 & 198 & 171 \\ \hline
  $r$ & 457 & 463 & 469 & 481 & 487 & 499 \\*[-1mm]
  $k$ & 133 & 21 & 37, 163 & 100, 211 & 232 & 139 \\ \hline
  $r$ & 511 & 523 & 541 & 547 & 553 & 559 \\*[-1mm]
  $k$ & 81, 137 & 60 & 129 & 40 & 23, 102 & 165, 178 \\ \hline
  $r$ & 571 & 577 & 589 & 601 & 607 & 613 \\*[-1mm]
  $k$ & 109 & 213 & 87, 273 & 24 & 210 & 65 \\ \hline
  $r$ & 619 & 631 & 637 & 643 & 661 & \\*[-1mm]
  $k$ & 252 & 43 & 165, 263 & 177 & 296 & \\ \hline
\end{tabular}
\end{center}
\captionsetup{width=11.2cm}
\caption{Continuation of table~\ref{rmn}:
  other values of $r$ and $k$ that produce groups $C_{r,1}^{(k)}$
  with order $3 r < 2\,000$.
  (For the values of $r$ in this table,
  only $l = 1$ produces group orders smaller than 2\,000.)
\label{rmn2}}
\end{table}
There is a very large number of groups $C_{n,l}^{(k)}$
of order smaller than 2\,000,
therefore we opt for giving their {\tt SmallGroups} identifiers
only in the appendix.

The generators of $C_{n,l}^{(k)}$ are the matrices $E$ in equation~\eqref{E},
$B_{n,k}$ in equation~\eqref{Bmk},
and $G_{n,r}$ in equation~\eqref{G},
where $r = n / l$.\footnote{For almost all the groups $C_{n,l}^{(k)}$
  of order smaller than 2\,000,
  the third generator $G_{n,r}$ is not really needed,
  \textit{i.e.}\ one may generate the group by using solely $E$
  and $B_{n,k}$.}

The groups $C_{n,l}^{(k)}$ only have singlet and triplet irreps.
The number of inequivalent singlet irreps is
three when $l$ cannot be divided by three and
nine when $l$ is a multiple of three.

\subsection{The groups $D_{3l,l}^{(1)}$}

We continue to use the notation of ref.~\cite{grimus}.
For an integer $l$ that is a multiple of three, 
the groups $D_{3l,l}^{(1)}$ have structure
$\left[ \left( \mathbbm{Z}_{3l} \times \mathbbm{Z}_l \right)
  \rtimes \mathbbm{Z}_3 \right] \rtimes \mathbbm{Z}_2$
and order $18 l^2$.
They are generated by the matrices $E$,
$I$,
and $B_{3l,1} = \mathrm{diag} \left( \nu,\, \nu,\, \nu^{-2} \right)$
for $\nu = \exp{\left[ 2 i \pi \!
    \left/ \left( 3l \right) \right. \right]}$.
There are only three groups $D_{3l,l}^{(1)}$ of order smaller than 2\,000:
\bs
\ba
D_{9,3}^{(1)} &\cong& \left[ 162, 14 \right], \\
D_{18,6}^{(1)} &\cong& \left[ 648, 259 \right], \\
D_{27,9}^{(1)} &\cong& \left[ 1458, 659 \right].
\ea
\es

The groups $D_{3l,l}^{(1)}$ have six inequivalent singlets
and three inequivalent doublets for any value of $l$.
Besides,
they have $6(l-1)$ inequivalent triplet irreps
and $l(l-3)/2 + 1$ inequivalent six-plets.

\subsection{The exceptional subgroups of $SU(3)$}

The groups $\Delta \left( 3 n^2 \right)$ and $C_{n,l}^{(k)}$
form
the class~C of finite subgroups of $SU(3)$.
The groups $\Delta \left( 6 n^2 \right)$ and $D_{3l,l}^{(1)}$
form
the class~D of finite subgroups of $SU(3)$.
Both classes~C and~D contain infinite numbers of subgroups.
Besides these infinite classes of subgroups,
$SU(3)$ has six `exceptional' subgroups;\footnote{The groups
  $\Sigma \left( 36 \times 1 \right)$,
  $\Sigma \left( 72 \times 1 \right)$,
  $\Sigma \left( 216 \times 1 \right)$,
  and $\Sigma \left( 360 \times 1 \right)$ are subgroups of $PSU(3)$,
  \textit{i.e.} of $SU(3)$ divided by its $\mathbbm{Z}_3$ center.
  They are \emph{not}\/ subgroups of $SU(3)$.}$^,$\footnote{The group
  $\Sigma \left( 60 \right)$ is in fact a subgroup of $SO(3)$,
  \textit{i.e.}\ it may be represented through real $3 \times 3$ matrices.}
their {\tt SmallGroups} identifiers are given in table~\ref{except}.
\begin{table}
\begin{center}
\renewcommand{\arraystretch}{1.2}
\begin{tabular}{|c|c|c|c|c|c|} \hline
  $\Sigma \left( 60 \right)$ &
  $\Sigma \left( 36\times3 \right)$ &
  $\Sigma \left( 168 \right)$ &
  $\Sigma \left( 72\times3 \right)$ &
  $\Sigma \left( 216\times3 \right)$ &
  $\Sigma \left( 360\times3 \right)$ \\*[-1mm]
  $\left[ 60, 5 \right]$ &
  $\left[ 108, 15 \right]$ &
  $\left[168, 42 \right]$ &
  $\left[ 216, 88 \right]$ &
  $\left[ 648, 532 \right]$ &
  $\left[ 1080, 260 \right]$ \\ \hline
\end{tabular}
\end{center}
\captionsetup{width=12.7cm}
\caption{The {\tt SmallGroups} identifiers
  of the exceptional subgroups of $SU(3)$.
\label{except}}
\end{table}
The generators of the exceptional subgroups are given,
for instance,
in ref.~\cite{ludl},
together with the references to the original papers.

The group $\Sigma \left( 60 \right)$ is isomorphic to $A_5$,
the group of the even permutations of five objects,
and to the symmetry group of the regular icosahedron and regular dodecahedron.
The group $\Sigma \left( 168 \right)$
is isomorphic to the projective special linear group $PSL \left( 2, 7 \right)$
and also to the general linear group $GL \left( 3, 2 \right)$.

The number of inequivalent $p$-dimensional irreps
of the exceptional finite subgroups of $SU(3)$
is given in table~\ref{exceptnumbers}~\cite{principal}.
\begin{table}
\begin{center}
\renewcommand{\arraystretch}{1.2}
\begin{tabular}{|c||c|c|c|c|c|c|} \hline
  group & $p=1$ & $p=2$ & $p=3$ & $p=4$ & $p=5$ & $p=6$
  \\ \hline
  $\Sigma \left( 60 \right)$ & 1 & 0 & 2 & 1 & 1 & 0 
  \\
  $\Sigma \left( 36\times3 \right)$ & 4 & 0 & 8 & 2 & 0 & 0
  \\
  $\Sigma \left( 168 \right)$ & 1 & 0 & 2 & 0 & 0 & 1
  \\
  $\Sigma \left( 72\times3 \right)$ & 4 & 1 & 8 & 0 & 0 & 2
  \\
  $\Sigma \left( 216\times3 \right)$ & 3 & 3 & 7 & 0 & 0 & 6
  \\
  $\Sigma \left( 360\times3 \right)$ & 1 & 0 & 4 & 0 & 2 & 2
  \\ \hline \hline
  group & $p=7$ & $p=8$ & $p=9$ & $p=10$ & $p=15$ &
  \\ \hline
  $\Sigma \left( 60 \right)$ & 0 & 0 & 0 & 0 & 0 & 
  \\
  $\Sigma \left( 36\times3 \right)$ & 0 & 0 & 0 & 0 & 0 &
  \\
  $\Sigma \left( 168 \right)$ & 1 & 1 & 0 & 0 & 0 &
  \\
  $\Sigma \left( 72\times3 \right)$ & 0 & 1 & 0 & 0 & 0 &
  \\
  $\Sigma \left( 216\times3 \right)$ & 0 & 3 & 2 & 0 & 0 &
  \\
  $\Sigma \left( 360\times3 \right)$ & 0 & 2 & 3 & 1 & 2 &
  \\ \hline
\end{tabular}
\end{center}
\captionsetup{width=11.4cm}
\caption{The number of inequivalent $p$-dimensional irreducible representations
  of the exceptional subgroups of $SU(3)$.
\label{exceptnumbers}}
\end{table}

\section{Finite subgroups of $U(3)$} \label{u3}

In this section we give
the generators and the {\tt SmallGroups} identifiers
of all the finite subgroups of $U(3)$ that
\begin{itemize}
  \item are not subgroups of $SU(3)$,
\item have a faithful three-dimensional irrep,
\item cannot be written as the direct product of a smaller group
  and a cyclic group,
\item have less than 2\,000 elements.
\end{itemize}
%
For most groups,
we also give the numbers of inequivalent irreps of each dimension.

There is at present no mathematical classification
of the finite subgroups of $U(3)$.
Therefore,
we will just classify the various subgroups that we have found
using the {\tt SmallGroups} library and {\tt GAP},
by constructing `series' of subgroups that have
generators, structures, and numbers of irreps
related among themselves.
Unfortunately,
there is some degree of ambiguity in this task,
since any group may always be generated by different sets of generators.
It is moreover often found that groups with related generators
end up having quite different structures.
Still,
we hope
to be able
to shed some light on the possible types of subgroups of $U(3)$.

\subsection{The generators} \label{U3sec1}

We firstly define some $3 \times 3$ matrices that often appear
as generators of the $U(3)$ subgroups.

Let
\begin{itemize}
\item $r$ be a product of prime numbers $p_1, p_2, \ldots$
  which are of the form $p_j = 6 i_j + 1$,
  where the numbers $i_j$ are integers;
\item $k$ be an integer which is a function of $r$ defined by
  $1 + k + k^2 = 0\ \mbox{mod}\ r$
  and $k \le \left. \left( r - 1 \right) \right/ 2$.
  For most values of $r$ there is only one possible $k$,
  but for some $r$ more than one $k$ are possible.
\end{itemize}
The lowest $r$ and the corresponding $k$ are given in table~\ref{tabrk}.
\begin{table}
\begin{center}
\renewcommand{\arraystretch}{1.2}
\begin{tabular}{|c||c|c|c|c|c|c|c|c|c|} \hline
  $r$ & 7 & 13 & 19 & 31 & 37 & 43 & 49 & 61 & 67 \\*[-1mm]
  $k$ & 2 & 3 & 7 & 5 & 10 & 6 & 18 & 13 & 29 \\
  \hline
  $r$ & 73 & 79 & 91 & 97 & 103 & 109 & 127 & 133 & 139 \\*[-1mm]
  $k$ & 8 & 23 & 9,\ 16 & 35 & 46 & 45 & 19 & 11,\ 30 & 42
  \\ \hline
  $r$ & 151 & 157 & 163 & 169 & 181 & 193 & 199 & 211 & 217 \\*[-1mm]
  $k$ & 32 & 12 & 58 & 22 & 48 & 84 & 92 & 14 & 25,\ 67 \\ \hline
\end{tabular}
\end{center}
\captionsetup{width=11.2cm}
\caption{The lowest possible values of $r$ and the corresponding values of $k$.
\label{tabrk}}
\end{table}
In this section,
whenever we let $r$ and $k$ denote a pair of integers,
we will be referring to one of the pairs in table~\ref{tabrk}.
The matrix
\be
\label{rho}
B_{r,k} = \mathrm{diag} \left( \rho,\, \rho^k,\, \rho^{-1-k} \right),
\quad \mbox{where}\ \rho = \exp{\left( 2 i \pi / r \right)},
\ee
appears as generator of many $U(3)$ subgroups.
Notice that $B_{r,k} \in SU(3)$.

We use the definition of $L_n$ in equation~\eqref{nu}.
Notice that $L_n \in SU(3)$.
The matrix
\be
\label{a4}
L_2 = \mathrm{diag} \left( 1,\, -1,\, -1 \right)
\ee
is especially useful.
We will also encounter
\be
L_3 = \mathrm{diag} \left( 1,\, \omega,\, \omega^2 \right),
\quad \mbox{where}\ \omega = \exp{\left( 2 i \pi / 3 \right)}.
\label{Pgen2}
\ee

Let $m$ be an integer.
We define
\bs
\ba
\label{mu}
E_m &\equiv&
\left( \begin{array}{ccc} 0 & \mu & 0 \\ 0 & 0 & \mu \\ \mu & 0 & 0
\end{array} \right),
\\
\label{Pgen3}
Z_m &\equiv& \left( \begin{array}{ccc} 0 & 0 & \mu \\
  1 & 0 & 0 \\ 0 & 1 & 0 \end{array}
\right), \\
T_1(m) &\equiv& \mbox{diag} \left( 1,\, \mu,\, \mu^2 \right),
\\
T_2(m) &\equiv& \mbox{diag} \left( 1,\, \mu^2,\, \mu \right),
\quad \mbox{where}\
\mu = \exp{\left[ 2 i \pi \left/\, \left( 3^m \right) \right. \right]}.
\ea
\es
The matrix $E \equiv E_0$ in equation~\eqref{E} is especially useful.
Both $E_0$ and $E_1$ have unit determinant,
but $E_m \notin SU(3)$ for $m>1$.
 
Let $m$ and $j$ be integers.
We define
\be
\label{xi}
F_{m,j} \equiv
\left( \begin{array}{ccc} 0 & 0 & - \xi \\ 0 & - \xi & 0 \\
  - \xi & 0 & 0 \end{array} \right),
\quad \mbox{where}\
\xi = \exp{\left[ 2 i \pi \left/\, \left( 3^m 2^j \right) \right. \right]}.
\ee
Notice that $F_{m,j} \notin SU(3)$ for $m \ge 2$ or $j \ge 1$.
The matrix $I \equiv F_{0,0}$ in equation~\eqref{Igen} has already been useful;
also useful is
\be
\label{Igen2}
I' \equiv F_{0,1} =
\left( \begin{array}{ccc} 0 & 0 & 1 \\ 0 & 1 & 0 \\ 1 & 0 & 0
  \end{array} \right) = -I.
\ee

Let $\omega = \exp{\left( 2 i \pi / 3 \right)}$
and $\mu = \exp{\left[ 2 i \pi \left/\, \left( 3^m \right) \right. \right]}$.
We define
\bs
\ba
X_1(m) &\equiv& \mathrm{diag} \left( \mu \omega,\, \mu \omega,\,
\mu \omega^2 \right),
\label{X1}
\\
X_2(m) &\equiv& \mathrm{diag} \left( \mu \omega^2,\, \mu \omega^2,\,
\mu \omega \right),
\label{X2}
\\
Y_1(m) &\equiv& \mathrm{diag} \left( \mu,\, \mu \omega,\,
\mu \omega^2 \right),
\label{y1} \\
Y_2(m) &\equiv& \mathrm{diag} \left( \mu,\, \mu \omega^2,\,
\mu \omega \right),
\label{y2} \\
X_3(m) = Y_3(m) &\equiv& \mathrm{diag} \left( \mu,\, \mu,\, \mu \right).
\label{X3}
\ea
\es

Let $\omega = \exp{\left( 2 i \pi / 3 \right)}$.
We define
\be
\label{genK}
K \equiv
\frac{- i}{\sqrt{3}} \left( \begin{array}{ccc}
  1 & 1 & 1 \\ 1 & \omega & \omega^2 \\
  1 & \omega^2 & \omega \end{array} \right).
\ee
Notice that $K \in SU(3)$.
Let furthermore
$\xi = \exp{\left[ 2 i \pi \left/\, \left( 3^m 2^j \right) \right. \right]}$.
We define
\be
\label{genQ}
Q_{m,j} \equiv
\frac{- i \xi}{\sqrt{3}} \left( \begin{array}{ccc}
  1 & \omega^2 & \omega^2 \\ \omega^2 & \omega^2 & 1 \\
  1 & \omega & 1 \end{array} \right).
\ee
Notice that $\det{Q_{m,j}} = \xi^3 \neq 1$ in general.

\subsection{The series of groups that Ludl has discovered}
\label{U3sec2}

Ludl~\cite{ludl} has proved the existence of the following series
of finite subgroups of $U(3)$.

\paragraph{Groups $T_r^{(k)}(m)$:}
The group $T_r^{(k)}(m)$,
where $m$ is an integer larger than 1,\footnote{If $m=1$,
  then $T_r^{(k)}(1) \cong C_{r,1}^{(k)}$ is a subgroup of $SU(3)$.}
has structure $\mathbbm{Z}_r \rtimes \mathbbm{Z}_{3^m}$ and order $3^m r$.
The groups $T_r^{(k)}(m)$ of order smaller than 2\,000
are given in table~\ref{trk}.
\begin{table}
\begin{center}
\renewcommand{\arraystretch}{1.4}
\begin{tabular}{|c||c|c|c|c|c|} \hline
  $T_r^{(k)}(m)$ &
  $T_7^{(2)}(2)$ &
  $T_7^{(2)}(3)$ &
  $T_7^{(2)}(4)$ &
  $T_7^{(2)}(5)$ &
  $T_{13}^{(3)}(2)$ \\*[-1mm]
identifier &
$\left[ 63, 1 \right]$ &
$\left[ 189, 1 \right]$ &
$\left[ 567, 1 \right]$ &
$\left[ 1701, 68 \right]$ &
$\left[ 117, 1 \right]$ \\ \hline
  $T_r^{(k)}(m)$ &
  $T_{13}^{(3)}(3)$ &
$T_{13}^{(3)}(4)$ &
$T_{19}^{(7)}(2)$ &
$T_{19}^{(7)}(3)$ &
$T_{19}^{(7)}(4)$ \\*[-1mm]
identifier &
$\left[ 351, 1 \right]$ &
$\left[ 1053, 16 \right]$ &
$\left[ 171, 1 \right]$ &
$\left[ 513, 1 \right]$ &
$\left[ 1539, 16 \right]$ \\ \hline
$T_r^{(k)}(m)$ &
$T_{31}^{(5)}(2)$ &
$T_{31}^{(5)}(3)$ &
$T_{37}^{(10)}(2)$ &
$T_{37}^{(10)}(3)$ &
$T_{43}^{(6)}(2)$ \\*[-1mm]
identifier &
$\left[ 279, 1 \right]$ &
$\left[ 837, 1 \right]$ &
$\left[ 333, 1 \right]$ &
$\left[ 999, 1 \right]$ &
$\left[ 387, 1 \right]$ \\ \hline
$T_r^{(k)}(m)$ &
$T_{43}^{(6)}(3)$ &
$T_{49}^{(18)}(2)$ &
$T_{49}^{(18)}(3)$ &
$T_{61}^{(13)}(2)$ &
$T_{61}^{(13)}(3)$ \\*[-1mm]
identifier &
$\left[ 1161, 6 \right]$ &
$\left[ 441, 1 \right]$ &
$\left[ 1323, 1 \right]$ &
$\left[ 549, 1 \right]$ &
$\left[ 1647, 6 \right]$ \\ \hline
$T_r^{(k)}(m)$ &
$T_{67}^{(29)}(2)$ &
$T_{67}^{(29)}(3)$ &
$T_{73}^{(8)}(2)$ &
$T_{73}^{(8)}(3)$ &
$T_{79}^{(23)}(2)$ \\*[-1mm]
identifier &
$\left[ 603, 1 \right]$ &
$\left[ 1809, 6 \right]$ &
$\left[ 657, 1 \right]$ &
$\left[ 1971, 6 \right]$ &
$\left[ 711, 1 \right]$ \\ \hline
$T_r^{(k)}(m)$ &
$T_{91}^{(9)}(2)$ &
$T_{91}^{(16)}(2)$ &
$T_{97}^{(35)}(2)$ &
$T_{103}^{(46)}(2)$ &
$T_{109}^{(45)}(2)$ \\*[-1mm]
identifier &
$\left[ 819, 4 \right]$ &
$\left[ 819, 3 \right]$ &
$\left[ 873, 1 \right]$ &
$\left[ 927, 1 \right]$ &
$\left[ 981, 1 \right]$ \\ \hline
$T_r^{(k)}(m)$ &
$T_{127}^{(19)}(2)$ &
$T_{133}^{(11)}(2)$ &
$T_{133}^{(30)}(2)$ &
$T_{139}^{(42)}(2)$ &
$T_{151}^{(32)}(2)$ \\*[-1mm]
identifier &
$\left[ 1143, 1 \right]$ &
$\left[ 1197, 3 \right]$ &
$\left[ 1197, 4 \right]$ &
$\left[ 1251, 1 \right]$ &
$\left[ 1359, 1 \right]$ \\ \hline
$T_r^{(k)}(m)$ &
$T_{157}^{(12)}(2)$ &
$T_{163}^{(58)}(2)$ &
$T_{169}^{(22)}(2)$ &
$T_{181}^{(48)}(2)$ &
$T_{193}^{(84)}(2)$ \\*[-1mm]
identifier &
$\left[ 1413, 1 \right]$ &
$\left[ 1467, 1 \right]$ &
$\left[ 1521, 1 \right]$ &
$\left[ 1629, 1 \right]$ &
$\left[ 1737, 1 \right]$ \\ \hline
$T_r^{(k)}(m)$ &
$T_{199}^{(92)}(2)$ &
$T_{211}^{(14)}(2)$ &
$T_{217}^{(25)}(2)$ &
$T_{217}^{(67)}(2)$ &
 \\*[-1mm]
identifier &
$\left[ 1791, 1 \right]$ &
$\left[ 1899, 1 \right]$ &
$\left[ 1953, 3 \right]$ &
$\left[ 1953, 4 \right]$ &
 \\ \hline
\end{tabular}
\end{center}
\captionsetup{width=12.0cm}
\caption{The {\tt SmallGroups} identifiers
  of the groups $T_r^{(k)}(m)$ with order smaller than 2\,000.
\label{trk}}
\end{table}
Each of these groups has two generators,
which may be chosen to be $B_{r,k}$ in equation~\eqref{rho}
and $E_m$ in equation~\eqref{mu}.

The groups $T_r^{(k)}(m)$ have $3^m$ inequivalent singlet irreps;
all the remaining irreps of those groups are triplets.

\paragraph{Groups $\Delta \left( 3 n^2, m \right)$:} The group
$\Delta \left( 3 n^2, m \right)$,
where the integer $n$ \emph{cannot be divided by 3}\/
and $m>1$,\footnote{If $m=1$,
  then $\Delta \left( 3 n^2, 1 \right) \cong \Delta \left( 3 n^2 \right)$
  is a subgroup of $SU(3)$.}
has structure $\left( \mathbbm{Z}_n \times \mathbbm{Z}_n \right)
\rtimes \mathbbm{Z}_{3^m}$ and order $3^m n^2$.
The groups $\Delta \left( 3 n^2, m \right)$
of order less than 2\,000 are listed in table~\ref{deltanm}.
\begin{table}
\begin{center}
\renewcommand{\arraystretch}{1.2}
\begin{tabular}{|c||c|c|c|c|c|} \hline
  $n, m$ &
  $2, 2$ &
  $2, 3$ &
  $2, 4$ &
  $2, 5$ &
  \\*[-1mm]
  identifier &
  $\left[ 36, 3 \right]$ &
  $\left[ 108, 3 \right]$ &
  $\left[ 324, 3 \right]$ &
  $\left[ 972, 3 \right]$ &
  \\ \hline
  $n, m$ &
  $4, 2$ &
  $4, 3$ &
  $4, 4$ &
  $5, 2$ &
  $5, 3$ \\*[-1mm]
  identifier &
  $\left[ 144, 3 \right]$ &
  $\left[ 432, 3 \right]$ &
  $\left[ 1296, 3 \right]$ &
  $\left[ 225, 3 \right]$ &
  $\left[ 675, 5 \right]$ \\ \hline
  $n, m$ &
  $7, 2$ &
  $7, 3$ &
  $8, 2$ &
  $8, 3$ &
  $10, 2$ \\*[-1mm]
  identifier &
  $\left[ 441, 7 \right]$ &
  $\left[ 1323, 14 \right]$ &
  $\left[ 576, 3 \right]$ &
  $\left[ 1728, 3 \right]$ &
  $\left[ 900, 66 \right]$ \\ \hline
  $n, m$ &
  $11, 2$ &
  $13, 2$ &
  $14, 2$ & & \\*[-1mm]
identifier &
$\left[ 1089, 3 \right]$ &
$\left[ 1521, 7 \right]$ &
$\left[ 1764, 91 \right]$ & & \\ \hline
\end{tabular}
\end{center}
\captionsetup{width=11.9cm}
\caption{The {\tt SmallGroups} identifiers
  of the groups $\Delta \left( 3 n^2, m \right)$ with order smaller than 2\,000.
\label{deltanm}}
\end{table}
The group $\Delta \left( 3 n^2, m \right)$ is generated
by the matrices $L_n$ in equation~\eqref{nu}
and $E_m$ in equation~\eqref{mu}.

The groups $\Delta \left( 3 n^2, m \right)$
have $3^m$ inequivalent singlet irreps;
all the remaining irreps of these groups are triplets.

\paragraph{Groups $S_4(j)$:} The group $S_4(j)$,
where $j>1$,\footnote{The group $S_4(1)
  \cong \Delta \left( 6\times2^2 \right)$
  is a subgroup of $SU(3)$.}
has structure $A_4 \rtimes \mathbbm{Z}_{2^j}$
and order $3\times2^{j+2}$.\footnote{The group $A_4$ has structure
  $\left( \mathbbm{Z}_2 \times \mathbbm{Z}_ 2 \right) \rtimes \mathbbm{Z}_3$.}
There are six groups $S_4(j)$ of order smaller than 2\,000;
they are given in table~\ref{S4m}.
\begin{table}
\begin{center}
\renewcommand{\arraystretch}{1.2}
\begin{tabular}{|c||c|c|c|} \hline
$j$ & 2 & 3 & 4 \\*[-1mm]
identifier &
$\left[ 48, 30 \right]$ &
$\left[ 96, 65 \right]$ &
$\left[ 192, 186 \right]$ \\ \hline
$j$ & 5 & 6 & 7 \\*[-1mm]
identifier  &
$\left[ 384, 581 \right]$ &
$\left[ 768, 1085351 \right]$ &
$\left[ 1536, 408544687 \right]$ \\ \hline
\end{tabular}
\end{center}
\captionsetup{width=10.6cm}
\caption{The {\tt SmallGroups} identifiers
  of the groups $S_4(j)$
  with order smaller than 2\,000.
\label{S4m}}
\end{table}
The group $S_4(j)$ is generated by the matrices
$E$ in equation~\eqref{E},
$L_2$ in equation~\eqref{a4},
and $- F_{0,j}$,
where $F_{m,j}$ is given in equation~\eqref{xi}.

The group $S_4(j)$ has $2^j$ inequivalent singlet irreps,
$2^{j-1}$ inequivalent doublet irreps,
$2^j$ inequivalent triplet irreps,
and no other irreps.

\paragraph{Groups $\Delta \left( 6 n^2, j \right)$:}
The group $\Delta \left( 6 n^2, j \right)$,
where\footnote{In Ludl's paper the existence of $\Delta \left( 6 n^2, j \right)$
  has been proved for integers $n$ not divisible by 3.
  We have verified,
  though,
  that $\Delta \left( 6 n^2, j \right)$ exists for every $n > 1$,
  at least when $3\times2^j n^2 < 2\,000$.} $n>1$
and $j>1$,\footnote{The groups $\Delta \left( 6 n^2, 1 \right)
  \cong \Delta \left( 6 n^2 \right)$ are subgroups of $SU(3)$.}
has structure $\left[ \left( \mathbbm{Z}_n \times \mathbbm{Z}_n \right)
  \rtimes \mathbbm{Z}_3 \right] \rtimes \mathbbm{Z}_{2^j}$
and order $3\times2^j n^2$.
But,
for $n=2$,
$\Delta \left( 6 n^2, j \right) = S_4(j)$,
therefore we only need to take into account $n \ge 3$;
there are then the 24 groups $\Delta \left( 6 n^2, j \right)$
with order less than 2\,000 given in table~\ref{Delta}.
\begin{table}
\begin{center}
\renewcommand{\arraystretch}{1.2}
\begin{tabular}{|c||c|c|c|} \hline
  $n, j$ & 3, 2 & 3, 3 & 3, 4 \\*[-1mm]
  identifier &
  $\left[ 108, 11 \right]$ &
  $\left[ 216, 17 \right]$ &
  $\left[ 432, 33 \right]$ \\ \hline
  $n, j$ & 3, 5 & 3, 6 & 4, 2 \\*[-1mm]
  identifier &
  $\left[ 864, 69 \right]$ &
  $\left[ 1728, 185 \right]$ &
  $\left[ 192, 182 \right]$ \\ \hline
  $n, j$ & 4, 3 & 4, 4 & 4, 5 \\*[-1mm]
  identifier &
  $\left[ 384, 571 \right]$ &
  $\left[ 768, 1085333 \right]$ &
  $\left[ 1536, 408544678 \right]$ \\ \hline
  $n, j$ & 5, 2 & 5, 3 & 5, 4 \\*[-1mm]
  identifier &
  $\left[ 300, 13 \right]$ &
  $\left[ 600, 45 \right]$ &
  $\left[ 1200, 183 \right]$ \\ \hline
  $n, j$ & 6, 2 & 6, 3 & 6, 4 \\*[-1mm]
  identifier &
  $\left[ 432, 260 \right]$ &
  $\left[ 864, 703 \right]$ &
  $\left[ 1728, 2855 \right]$ \\ \hline
  $n, j$ & 7, 2 & 7, 3 & 8, 2 \\*[-1mm]
  identifier &
  $\left[ 588, 16 \right]$ &
  $\left[ 1176, 57 \right]$ &
  $\left[ 768, 1085335 \right]$ \\ \hline
  $n, j$ & 8, 3 & 9, 2 & 9, 3 \\*[-1mm]
  identifier &
  $\left[ 1536, 408544641 \right]$ &
  $\left[ 972, 64 \right]$ &
  $\left[ 1944, 70 \right]$ \\ \hline
  $n, j$ & 10, 2 & 11, 2 & 12, 2 \\*[-1mm]
  identifier &
  $\left[ 1200, 682 \right]$ &
  $\left[ 1452, 11 \right]$ &
  $\left[ 1728, 2847 \right]$ \\ \hline
\end{tabular}
\end{center}
\captionsetup{width=12.0cm}
\caption{The {\tt SmallGroups} identifiers
  of the groups $\Delta \left( 6n^2, j \right)$
  with order smaller than 2\,000.
\label{Delta}}
\end{table}
The generators of $\Delta \left( 6 n^2, j \right)$ are the matrices
$E$ in equation~\eqref{E},
$L_n$ in equation~\eqref{nu},
and $- F_{0,j}$.
It is clear that $\Delta \left( 6 n^2, j \right)$
is just a generalization of $S_4(j)$ for $n > 2$.

The groups $\Delta \left( 6 n^2, j \right)$
have $2^j$ inequivalent singlet irreps
and $2^j (n-1)$ inequivalent triplet irreps
for any value of $n$.
When $n$ cannot be divided by three,
those groups have,
besides,
$2^{j-1}$ doublet irreps;
when $n$ is a multiple of three,
the number of inequivalent doublet irreps is $2^{j+1}$.
All the remaining irreps of those groups are six-plets.

\paragraph{Groups $\Delta' \left( 6 n^2, m, j \right)$:}
These groups,
where $n$ can be divided by 3,
$m \ge 2$,\footnote{The groups $\Delta' \left( 6 n^2, 1, j \right)$
are the same as the groups $\Delta \left( 6 n^2, j \right)$.}
and $j \ge 1$,
have structure\footnote{The exception is $\Delta' \left(
  6\times9^2, 2, 1 \right)$,
  which has structure $\left[ \left( \mathbbm{Z}_9 \times \mathbbm{Z}_9
    \times \mathbbm{Z}_3 \right) \rtimes \mathbbm{Z}_3 \right]
  \rtimes \mathbbm{Z}_2$ instead of
  $\left[ \left( \mathbbm{Z}_{27} \times \mathbbm{Z}_9
    \right) \rtimes \mathbbm{Z}_3 \right] \rtimes \mathbbm{Z}_2$.}
  $\left[ \left(
  \mathbbm{Z}_{3^{m-1} n} \times \mathbbm{Z}_n \right)
  \rtimes \mathbbm{Z}_3 \right] \rtimes \mathbbm{Z}_{2^j}$
and order $3^m 2^j n^2$.
There are the 12 groups with  order less than 2\,000 in table~\ref{Deltaprime}.
\begin{table}
\begin{center}
\renewcommand{\arraystretch}{1.2}
\begin{tabular}{|c||c|c|c|c|} \hline
  $n, m, j$ & 3, 2, 1 & 3, 2, 2 & 3, 2, 3 & 3, 2, 4 \\*[-1mm]
  identifier &
  $\left[ 162, 44 \right]$ &
  $\left[ 324, 102 \right]$ &
  $\left[ 648, 244 \right]$ &
  $\left[ 1296, 647 \right]$
  \\ \hline
  $n, m, j$ & 3, 3, 1 & 3, 3, 2 & 3, 3, 3 & 3, 4, 1 \\*[-1mm]
  identifier &
  $\left[ 486, 164 \right]$ &
  $\left[ 972, 348 \right]$ &
  $\left[ 1944, 746 \right]$ &
  $\left[ 1458, 1354 \right]$
  \\ \hline
  $n, m, j$ & 6, 2, 1 & 6, 2, 2 & 6, 3, 1 & 9, 2, 1 \\*[-1mm]
  identifier &
  $\left[ 648, 563 \right]$ &
  $\left[ 1296, 2113 \right]$ &
  $\left[ 1944, 2415 \right]$ &
  $\left[ 1458, 1371 \right]$
  \\ \hline
\end{tabular}
\end{center}
\captionsetup{width=11.6cm}
\caption{{\tt SmallGroups} identifiers
  of the groups $\Delta' \left( 6n^2, m, j \right)$
  with order smaller than 2\,000.
\label{Deltaprime}}
\end{table}
The generators of $\Delta' \left( 6 n^2, m, j \right)$ are the matrices
$E$ in equation~\eqref{E},
$L_n$ in equation~\eqref{nu},
and $- F_{m,j}$ in equation~\eqref{xi}.

The groups $\Delta' \left( 6 n^2, m, j \right)$ have
$3^{m-1}2^j$ inequivalent singlet irreps
and $3^{m-1}2^{j+1}$ inequivalent doublet irreps.
There are also $\left( n - 1 \right) 3^{m-1} 2^j$ inequivalent triplets;
the remaining irreps are six-plets.

\subsection{New series of groups that we have discovered}
\label{U3sec3}

Ludl~\cite{ludl} has derived the existence
of the series of groups in the previous subsection
by applying mathematical theorems that he demonstrated.
We have discovered some further series of groups
through a careful inspection
of the list of all the finite subgroups of $U(3)$ of order smaller than 2\,000
that we have produced,
together with some guesswork.
Clearly,
since there are no theorems supporting our method,
we cannot be sure that our series of groups
extend to groups of order larger than 2\,000.
Still,
the series of groups in this subsection seem to us to be on firm standing,
since they are quite large and display no exceptions up to group order 2\,000.

\paragraph{Groups $L_r^{(k)} \left( n, m \right)$:}
For an integer $n$ that cannot be divided by 3 and for $m > 1$,
these are groups with structure
$\left( \mathbbm{Z}_{rn} \times \mathbbm{Z}_n \right) \rtimes \mathbbm{Z}_{3^m}$
and order $3^m r n^2$.
While the groups $T_r^{(k)}(m)$ are generated by the matrices $B_{r,k}$ and $E_m$,
and the groups $\Delta \left( 3 n^2, m \right)$
are generated by the matrices $L_n$ and $E_m$,
the groups $L_r^{(k)} \left( n, m \right)$
are generated by all three matrices $B_{r,k}$,
$L_n$,
and $E_m$.
Thus,
the groups $L_r^{(k)} \left( n, m \right)$ simultaneously generalize
$T_r^{(k)}(m) = L_r^{(k)} \left( 1 , m \right)$
and $\Delta \left( 3 n^2, m \right) = L_1^{(0)} \left( n, m \right)$.
The groups $L_r^{(k)} \left( n, m \right)$ of order smaller than 2\,000
are listed in table~\ref{gamma}.
\begin{table}
\begin{center}
\renewcommand{\arraystretch}{1.4}
\begin{tabular}{|c||c|c|c|c|} \hline
  $L_r^{(k)} \left( n, m \right)$ &
  $L_7^{(2)} \left( 2, 2 \right)$ &
  $L_7^{(2)} \left( 2, 3 \right)$ &
  $L_7^{(2)} \left( 4, 2 \right)$ &
  $L_7^{(2)} \left( 5, 2 \right)$ \\*[-1mm]
identifier &
$\left[ 252, 11 \right]$ &
$\left[ 756, 11 \right]$ &
$\left[ 1008, 57 \right]$ &
$\left[ 1575, 7 \right]$ \\ \hline
  $L_r^{(k)} \left( n, m \right)$ &
  $L_{13}^{(3)} \left( 2, 2 \right)$ &
  $L_{13}^{(3)} \left( 2, 3 \right)$ &
  $L_{13}^{(3)} \left( 4, 2 \right)$ &
  $L_{19}^{(7)} \left( 2, 2 \right)$ \\*[-1mm]
identifier &
$\left[ 468, 14 \right]$ &
$\left[ 1404, 14 \right]$ &
$\left[ 1872, 60 \right]$ &
$\left[ 684, 11 \right]$ \\ \hline
  $L_r^{(k)} \left( n, m \right)$ &
  $L_{31}^{(5)} \left( 2, 2 \right)$ &
  $L_{37}^{(10)} \left( 2, 2 \right)$ &
  $L_{43}^{(6)} \left( 2, 2 \right)$ &
  $L_{49}^{(18)} \left( 2, 2 \right)$ \\*[-1mm]
identifier &
$\left[ 1116, 11 \right]$ &
$\left[ 1332, 14 \right]$ &
$\left[ 1548, 11 \right]$ &
$\left[ 1764, 11 \right]$ \\ \hline
\end{tabular}
\end{center}
\captionsetup{width=11.0cm}
\caption{{\tt SmallGroups} identifiers
  of the groups $L_r^{(k)} \left( n, m \right)$
  with order smaller than 2\,000.
\label{gamma}}
\end{table}

The groups $L_r^{(k)} \left( n, m \right)$ have $3^m$ inequivalent singlets;
the remaining irreps are triplets.

\paragraph{Groups $P_r^{(k)}(m)$,
  $Q_r^{(k)}(m)$,
  and $Q_r^{(k)\prime}(m)$:}
These groups exist for integer $m>1$ and have order $3^{m+1} r$.
The groups $P_r^{(k)} (m)$
have structure $\left( \mathbbm{Z}_ r \times \mathbbm{Z}_{3^m} \right)
\rtimes \mathbbm{Z}_3$;
the groups $Q_r^{(k)} (m)$
and $Q_r^{(k)\prime} (m)$
have structure $\mathbbm{Z}_{3^mr} \rtimes \mathbbm{Z}_3$.
The groups of order smaller than 2\,000
are listed in table~\ref{Gamma_rk}.
\begin{table}
\begin{center}
\renewcommand{\arraystretch}{1.4}
\begin{tabular}{|c||c|c|c|c|c|} \hline
  $_{\, r}^{(k)} (m)$ &
  $^{(2)}_{\, 7} (2)$ &
  $^{(3)}_{13} (2)$ &
  $^{(7)}_{19} (2)$ &
  $^{(5)}_{31} (2)$ &
  $^{(10)}_{37} (2)$  
\\
  $P^{(k)}_r (m)$ &
$\left[ 189, 7 \right]$ &
$\left[ 351, 7 \right]$ &
$\left[ 513, 8 \right]$ &
$\left[ 837, 7 \right]$ &
$\left[ 999, 8 \right]$
\\
  $Q^{(k)}_r (m)$ &
$\left[ 189, 4 \right]$ &
$\left[ 351, 4 \right]$ &
$\left[ 513, 5 \right]$ &
$\left[ 837, 5 \right]$ &
$\left[ 999, 6 \right]$
\\
  $Q^{(k) \prime}_r (m)$ &
$\left[ 189, 5 \right]$ &
$\left[ 351, 5 \right]$ &
$\left[ 513, 6 \right]$ &
$\left[ 837, 4 \right]$ &
$\left[ 999, 5 \right]$
\\ \hline
$^{(k)}_{\, r} (m)$ &
$^{(6)}_{43} (2)$ &
$^{(18)}_{49} (2)$ &
$^{(13)}_{61} (2)$ &
$^{(29)}_{67} (2)$ &
$^{(8)}_{73} (2)$
\\
$P^{(k)}_r (m)$ &
$\left[ 1161, 12 \right]$ &
$\left[ 1323, 7 \right]$ &
$\left[ 1647, 17 \right]$ &
$\left[ 1809, 17 \right]$ &
$\left[ 1971, 17 \right]$
\\
$Q^{(k)}_r (m)$ &
$\left[ 1161, 10 \right]$ &
$\left[ 1323, 4 \right]$ &
$\left[ 1647, 10 \right]$ &
$\left[ 1809, 10 \right]$ &
$\left[ 1971, 11 \right]$
\\
$Q^{(k) \prime}_r (m)$ &
$\left[ 1161, 11 \right]$ &
$\left[ 1323, 5 \right]$ &
$\left[ 1647, 11 \right]$ &
$\left[ 1809, 11 \right]$ &
$\left[ 1971, 10 \right]$
\\ \hline
  $^{(k)}_{\, r} (m)$ &
  $^{(2)}_{\, 7} (3)$ &
  $^{(3)}_{13} (3)$ &
  $^{(7)}_{19} (3)$ &
  $^{(2)}_{\, 7} (4)$ &  
\\
$P^{(k)}_r (m)$ &
$\left[ 567, 7 \right]$ &
$\left[ 1053, 27 \right]$ &
$\left[ 1539, 27 \right]$ &
$\left[ 1701, 128 \right]$ &
\\
$Q^{(k)}_r (m)$ &
$\left[ 567, 4 \right]$ &
$\left[ 1053, 26 \right]$ &
$\left[ 1539, 26 \right]$ &
$\left[ 1701, 127 \right]$ &
\\
$Q^{(k) \prime}_r (m)$ &
$\left[ 567, 5 \right]$ &
$\left[ 1053, 25 \right]$ &
$\left[ 1539, 25 \right]$ &
$\left[ 1701, 126 \right]$ &
\\
\hline
\end{tabular}
\end{center}
\captionsetup{width=12.6cm}
\caption{The {\tt SmallGroups} identifiers of the groups
  $P^{(k)}_r (m)$,
  $Q^{(k)}_r (m)$,
  and $Q^{(k)\prime}_r (m)$ with order smaller than 2\,000.
\label{Gamma_rk}}
\renewcommand{\arraystretch}{1.0}
\end{table}
The group $P^{(k)}_r (m)$
is generated by $B_{r,k}$ together with $L_3$ and $Z_{m-1}$.
The groups $Q^{(k)}_r (m)$ and $Q^{(k)\prime}_r (m)$
are generated by the matrices $B_{r,k}$ and $E$ together with
$Y_1(m)$ for $Q_r^{(k)} (m)$ or $Y_2(m)$ for $Q_r^{(k)\prime} (m)$.

The groups $P_r^{(k)}(m)$,
$Q_r^{(k)}(m)$,
and $Q_r^{(k)\prime}(m)$ have $3^m$ inequivalent singlets;
all their remaining irreps are triplets.

\paragraph{Groups $X (n)$:}

There are several groups
that have a three-dimensional irrep
where all the matrices are of one of the following types~\cite{ludl}:
\bs
\label{sutyo}
\ba
R \left( n, a, b, c \right) &\equiv&
\left( \begin{array}{ccc}  0 & 0 & \nu^a \\ \nu^b & 0 & 0 \\
  0 & \nu^c & 0 \end{array} \right), \\
V \left( n, a, b, c \right) &\equiv&
\left( \begin{array}{ccc}  0 & \nu^a & 0 \\ 0 & 0 & \nu^b \\
  \nu^c & 0 & 0 \end{array} \right), \\
W \left( n, a, b, c \right) &\equiv&
\left( \begin{array}{ccc}  \nu^a & 0 & 0 \\ 0 & \nu^b & 0 \\
  0 & 0 & \nu^c \end{array} \right),
\ea
\es
where $\nu = \exp{\left( 2 i \pi / n \right)}$.
We call them `groups RVW'.
The groups $X (n)$ are groups RVW where
\begin{itemize}
\item $n$ is a multiple of 3,
  \item the matrices $R \left( n, a, b, c \right)$
    have $a+ b + c = \left( n/ 3 \right) \mbox{mod}\ n$,
    \item the matrices $V \left( n, a, b, c \right)$
      have $a+ b + c = \left( 2 n / 3 \right) \mbox{mod}\ n$,
      \item the matrices $W \left( n, a, b, c \right)$
        have $a+ b + c = 0\ \mbox{mod}\ n$.
\end{itemize}
The groups $X (n)$ have order $3 n^2$;
the identifiers of the groups of order less than 2\,000 are in table~\ref{Xn}.
\begin{table}
\begin{center}
\renewcommand{\arraystretch}{1.2}
\begin{tabular}{|c||c|c|c|c|} \hline
$n$ & 3 & 6 & 9 & 12 \\*[-1mm]
identifier & $\left[ 27, 4 \right]$ & $\left[ 108, 21 \right]$ &
$\left[ 243, 27 \right]$ & $\left[ 432, 102 \right]$ \\ \hline
$n$ & 15 & 18 & 21 & 24 \\*[-1mm]
identifier & $\left[ 675, 11 \right]$ & $\left[ 972, 123 \right]$ &
$\left[ 1323, 42 \right]$ & $\left[ 1728, 1290 \right]$ \\ \hline  
\end{tabular}
\end{center}
\captionsetup{width=10.6cm}
\caption{The {\tt SmallGroups} identifiers of the groups $X (n)$
  with order smaller than 2\,000.
\label{Xn}}
\end{table}
The structure of $X (n)$ is
$\left[ \left( \mathbbm{Z}_{n/3} \times \mathbbm{Z}_{n/3} \right)
  \rtimes \mathbbm{Z}_9 \right] \rtimes \mathbbm{Z}_3$
provided $n$ is not a multiple of~9;
otherwise it is more complicated.
The groups $X (n)$ are generated by the matrices $L_n$ in equation~\eqref{nu}
and $Z_1$ in equation \eqref{Pgen3}.

The groups $X(n)$ have nine inequivalent singlets;
their remaining irreps are all triplets.

\subsection{Tentative series of groups} \label{U3sec4}

We have found a few more series of groups
through inspection of the list of the finite subgroups of $U(3)$
of order less than 2\,000.
However,
these series have few groups each and we can hardly ascertain whether and how
they extend to groups of order larger than 2\,000.

\paragraph{Groups $S_r^{(k)} (m)$, $S_r^{(k)\prime} (m)$, $Y_r^{(k)} (m)$,
  and $V_r^{(k)} (m)$:}
These groups exist for $m \ge 2$.
The groups $S_r^{(k)} (m)$ and $S_r^{(k)\prime} (m)$ have structure
$\left( \mathbbm{Z}_{3^{m} r} \times \mathbbm{Z}_{3}
\right) \rtimes \mathbbm{Z}_3$;
the groups $Y_r^{(k)} (m)$ have structure
$\left( \mathbbm{Z}_{3^{m-1} r} \times \mathbbm{Z}_{3}
\times \mathbbm{Z}_{3} \right) \rtimes \mathbbm{Z}_3$;
the groups $V_r^{(k)} (m)$ have structure
$\mathbbm{Z}_r \rtimes
\left[ \left( \mathbbm{Z}_{3^{m-1}} \times \mathbbm{Z}_{3} \right)
  . \left( \mathbbm{Z}_3 \times \mathbbm{Z}_{3} \right) \right]$;
they all have order $3^{m+2} r$.
The generators are the matrices $B_{r,k}$ in equation~\eqref{rho},
together with
\begin{itemize}
\item $E$ in equation~\eqref{E},
  $L_3$ in equation~\eqref{Pgen2},
  and $X_3(m)$ in equation~\eqref{X3} for $S_r^{(k)} (m)$;
\item $E$,
  $L_3$,
  and $X_1(m)$ in equation~\eqref{X1} for $S_r^{(k)\prime} (m)$;
\item $E$,
  $X_1 (m-1)$,
  and $X_3 (m-1)$ for $Y_r^{(k)} (m)$;
\item $X_2(2)$ in equation~\eqref{X2},
  $Z_1$ in equation~\eqref{Pgen3},
  and $L_{3^{m-1}}$
  for $V_r^{(k)} (m)$.
\end{itemize}
The groups $S_r^{(k)} (m)$, $S_r^{(k)\prime} (m)$, $Y_r^{(k)} (m)$,
and $V_r^{(k)} (m)$ of order less than 2\,000 are in table~\ref{g}.
\begin{table}
\begin{center}
\renewcommand{\arraystretch}{1.4}
\begin{tabular}{|c|c|c|c|} \hline
  $S_7^{(2)} (2)$ &
  $S_{13}^{(3)} (2)$ &
  $S_{19}^{(7)} (2)$ &
  $S_7^{(2)} (3)$ \\*[-1mm]
$\left[ 567, 36 \right]$ &
$\left[ 1053, 47 \right]$ &
$\left[ 1539, 47 \right]$ &
$\left[ 1701, 240 \right]$ \\ \hline
  $S_7^{(2)\prime} (2)$ &
  $S_{13}^{(3)\prime} (2)$ &
  $S_{19}^{(7)\prime} (2)$ &
  $S_7^{(2)\prime} (3)$ \\*[-1mm]
$\left[ 567, 12 \right]$ &
$\left[ 1053, 32 \right]$ &
$\left[ 1539, 32 \right]$ &
$\left[ 1701, 115 \right]$ \\ \hline
  $Y_7^{(2)} (2)$ &
  $Y_{13}^{(3)} (2)$ &
  $Y_{19}^{(7)} (2)$ &
  $Y_7^{(2)} (3)$ \\*[-1mm]
$\left[ 567, 23 \right]$ &
$\left[ 1053, 29 \right]$ &
$\left[ 1539, 29 \right]$ &
$\left[ 1701, 261 \right]$ \\ \hline
  $V_7^{(2)} (2)$ &
  $V_{13}^{(3)} (2)$ &
  $V_{19}^{(7)} (2)$ &
  $V_7^{(2)} (3)$ \\*[-1mm]
$\left[ 567, 14 \right]$ &
$\left[ 1053, 37 \right]$ &
$\left[ 1539, 37 \right]$ &
$\left[ 1701, 138 \right]$ \\ \hline
\end{tabular}
\end{center}
\captionsetup{width=8.2cm}
\caption{The {\tt SmallGroups} identifiers of the groups $S_r^{(k)} (m)$,
  $S_r^{(k)\prime} (m)$,
  $Y_r^{(k)} (m)$,
  and $V_r^{(k)} (m)$ with order smaller than 2\,000.
\label{g}}
\end{table}

The groups $S_r^{(k)}(m)$ have $3^{m+1}$ inequivalent singlets.
The groups $S_r^{(k)\prime}(m)$ and $Y_r^{(k)}(m)$ have $3^m$ inequivalent singlets.
The groups $V_r^{(k)}(m)$ have nine inequivalent singlets.
All the remaining irreps of all those groups are triplets.

\paragraph{Groups $M_r^{(k)}$,
  $M_r^{(k) \prime}$,
  and $J_r^{(k)}$:}
These groups have order $108 r$.
The groups $M_r^{(k)}$ and $M_r^{(k) \prime}$ have structure
$\left( \mathbbm{Z}_{18 r} \times \mathbbm{Z}_{2} \right) \rtimes \mathbbm{Z}_3$;
the groups $J_r^{(k)}$ have structure
$\left[ \left( \mathbbm{Z}_{2 r} \times \mathbbm{Z}_{2} \right)
\rtimes \mathbbm{Z}_9 \right] \rtimes \mathbbm{Z}_3$.
The generators are the matrices $B_{r,k}$ in equation~\eqref{rho}
and $L_2$ in equation~\eqref{a4} together with
\begin{itemize}
\item $E$ in equation~\eqref{E} and $Y_1(2)$ in equation~\eqref{y1}
  for $M_r^{(k)}$,
\item $E$ and $Y_2(2)$ in equation~\eqref{y2} for $M_r^{(k) \prime}$,
\item $L_3$ in equation~\eqref{Pgen2} and $Z_1$ in equation~\eqref{Pgen3}
  for $J_r^{(k)}$.
\end{itemize}
The groups $M_r^{(k)}$,
$M_r^{(k) \prime}$,
and $J_r^{(k)}$ of order less than 2\,000
are in table~\ref{Psi}.
\begin{table}
\begin{center}
\renewcommand{\arraystretch}{1.4}
\begin{tabular}{|c|c|c|c|c|c|} \hline
  $M_7^{(2)}$ &
  $M_7^{(2) \prime}$ &
  $J_7^{(2)}$ &
  $M_{13}^{(3)}$ &
  $M_{13}^{(3) \prime}$ &
  $J_{13}^{(3)}$ \\*[-1mm]
$\left[ 756, 113 \right]$ &
$\left[ 756, 114 \right]$ &
$\left[ 756, 116 \right]$ &
$\left[ 1404, 137 \right]$ &
$\left[ 1404, 138 \right]$ &
$\left[ 1404, 140 \right]$ \\ \hline
\end{tabular}
\end{center}
\captionsetup{width=13.2cm}
\caption{The {\tt SmallGroups} identifiers of the groups $M_r^{(k)}$,
  $M_r^{(k)\prime}$,
  and $J_r^{(k)}$
  with order smaller than 2\,000.
\label{Psi}}
\end{table}

Each of the groups $M_r^{(k)}$,
$M_r^{(k)\prime}$,
and $J_r^{(k)}$ has nine inequivalent singlets.
All the remaining irreps of those groups are triplets.

\paragraph{Groups $W \left( n, m \right)$:}
The groups $W \left( n, m \right)$,
where $n$ cannot be divided by 3 and $m>1$,
are generated by the matrices~$E$ in equation~\eqref{E},
$L_n$ in equation~\eqref{nu},
and $Y_1(m)$ in equation~\eqref{y1}.
They have structure $\left( \mathbbm{Z}_{3^m n} \times \mathbbm{Z}_{n} 
\right) \rtimes \mathbbm{Z}_3$
and order $3^{m+1} n^2$.
The groups $W \left( n, m \right)$
with order smaller than 2\,000 are listed in table~\ref{g1}.
\begin{table}
\begin{center}
\renewcommand{\arraystretch}{1.2}
\begin{tabular}{|c||c|c|c|c|} \hline
$n, m$ & 1, 2 & 1, 3 & 1, 4 & 1, 5 \\*[-1mm]
identifier & $\left[ 27, 4 \right]$ & $\left[ 81, 6 \right]$ &
$\left[ 243, 24 \right]$ & $\left[ 729, 94 \right]$ \\ \hline
$n, m$ & 2, 2 & 2, 3 & 2, 4 & 4, 2 \\*[-1mm]
identifier & $\left[ 108, 19 \right]$ & $\left[ 324, 43 \right]$ &
$\left[ 972, 117 \right]$ & $\left[ 432, 100 \right]$ \\ \hline
$n, m$ & 4, 3 & 5, 2 & 7, 2 & 8, 2 \\*[-1mm]
identifier & $\left[ 1296, 220 \right]$ & $\left[ 675, 9 \right]$ &
$\left[ 1323, 40 \right]$ & $\left[ 1728, 1286 \right]$ \\ \hline
\end{tabular}
\end{center}
\captionsetup{width=10.8cm}
\caption{{\tt SmallGroups} identifiers of the groups $W \left( n, m \right)$
  with order smaller than 2\,000.
\label{g1}}
\end{table}

Each of the groups $W \left( n, m \right)$ has $3^m$ inequivalent singlets;
the remaining irreps of those groups are triplets.

\paragraph{Groups $Z \left( n, m \right)$,
  $Z' \left( n, m \right)$,
  and $Z'' \left( n, m \right)$:}
These groups,
where $n$ is a multiple of 3 and $m>1$,
have structure\footnote{An exception are the groups
  $Z \left( 9, 2 \right)$ and $Z' \left( 9, 2 \right)$,
  which happen to be isomorphic and are of the form
  $\left( \mathbbm{Z}_9 \times \mathbbm{Z}_9 \times \mathbbm{Z}_3 \right)
  \rtimes \mathbbm{Z}_3$ instead of
  $\left( \mathbbm{Z}_{27} \times \mathbbm{Z}_9 \right)
  \rtimes \mathbbm{Z}_3$.}
$\left( \mathbbm{Z}_{3^{m-1} n} \times \mathbbm{Z}_{n} \right)
\rtimes \mathbbm{Z}_3$
and order $3^m n^2$.
The groups with order smaller than 2\,000
are listed in table~\ref{gI}.
\begin{table}
\begin{center}
\renewcommand{\arraystretch}{1.2}
\begin{tabular}{|c||c|c|c|c|} \hline
  $n, m$ &
  $3, 2$ &
  $6, 2$ &
  $9, 2$ &
  $12, 2$ \\*[-1mm]
  $Z \left( n, m \right)$ &
  $\left[ 81, 14 \right]$ &
  $\left[ 324, 128 \right]$ &
  $\left[ 729, 397 \right]$ &
  $\left[ 1296, 1499 \right]$ \\*[-1mm]
  $Z' \left( n, m \right)$ &
  $\left[ 81, 8 \right]$ &
  $\left[ 324, 49 \right]$ &
  $\left[ 729, 397 \right]$ &
  $\left[ 1296, 227 \right]$ \\*[-1mm]
  $Z'' \left( n, m \right)$ &
  $C_{9,3}^{(1)}$ &
  $C_{18,6}^{(1)}$ &
  $\Delta \left( 3\times9^2 \right)$ &
  $C_{36,12}^{(1)}$ \\ \hline
  $n, m$ &
  $3, 3$ &
  $6, 3$ &
  $3, 4$ & \\*[-1mm]
  $Z \left( n, m \right)$ &
  $\left[ 243, 50 \right]$ &
  $\left[ 972, 520 \right]$ &
  $\left[ 729, 393 \right]$ & \\*[-1mm]
  $Z' \left( n, m \right)$ &
  $\left[ 243, 20 \right]$ &
  $\left[ 972, 152 \right]$ &
  $\left[ 729, 64 \right]$ & \\*[-1mm]
  $Z'' \left( n, m \right)$ &
  $\left[ 243, 19 \right]$ &
  $\left[ 972, 153 \right]$ &
  $\left[ 729, 63 \right]$ & \\ \hline
\end{tabular}
\end{center}
\captionsetup{width=11.0cm}
\caption{{\tt SmallGroups} identifiers
  of the groups $Z \left( n, m \right)$,
  $Z' \left( n, m \right)$,
  and $Z'' \left( n, m \right)$
  with order smaller than 2\,000.
\label{gI}}
\end{table}
The generators of $Z \left( n, m \right)$ are
just the same as those of $W \left( n, m \right)$,
\textit{viz.}\ $E$,
$L_n$,
and $Y_1(m)$---the only
difference being that for $Z \left( n, m \right)$
the integer $n$ is a multiple of 3
while for $W \left( n, m \right)$
the integer $n$ cannot be divided by 3.
The groups $Z' \left( n, m \right)$
are generated by the matrices~$E$,
$L_n$,
and~$X_1(m)$.
The groups $Z'' \left( n, m \right)$
are generated by the matrices~$E$,
$L_n$,
and~$X_2(m)$.
Notice that,
for $m = 2$,
$Z'' \left( n, m \right)$ is generated by matrices
with unit determinant and therefore
it is a subgroup of $SU(3)$.

Each of the groups $Z \left( n, m \right)$
and $Z'' \left( n, m \right)$ has $3^{m+1}$ inequivalent singlets.
The groups $Z' \left( n, m \right)$ have $3^m$ inequivalent singlets.
All the remaining irreps of all those groups are triplets.

\paragraph{Groups $Z \left( n, m, j \right)$
  and $Z' \left(  n, m, j \right)$:}
These groups,
where
\begin{itemize}
\item $n$ is a multiple of 3,
\item $m>1$,
\item $j$ is an integer,
\end{itemize}
have order $3^m 2^j n^2$. 
The groups $Z \left( n, m, j \right)$ and $Z' \left( n, m, j \right)$
with order smaller than 2\,000
are in table~\ref{gg}.\footnote{The group $Z \left( 9, 2, 1 \right)$
  is isomorphic to the group $\Delta' \left( 6\times9^2, 2, 1 \right)$
  and we omit it from table~\ref{gg},
  since it
  has numbers of irreps quite inconsistent with the ones
  of the other groups $Z \left( n, m, j \right)$.}
\begin{table}
\begin{center}
\renewcommand{\arraystretch}{1.2}
\begin{tabular}{|c||c|c|c|c|} \hline
  $n, m, j$ &
  $3, 2, 1$ &
  $3, 2, 2$ &
  $3, 2, 3$ &
  $3, 2, 4$ \\*[-1mm]
  $Z \left( n, m, j \right)$  &
  $\left[ 162, 12 \right]$ &
  $\left[ 324, 15 \right]$ &
  $\left[ 648, 21 \right]$ &
  $\left[ 1296, 37 \right]$ \\*[-1mm]
  $Z' \left( n, m, j \right)$ &
  $D_{9,3}^{(1)}$ &
  $\left[ 324, 17 \right]$ &
  $\left[ 648, 23 \right]$ &
  $\left[ 1296, 39 \right]$ \\ \hline
  $n, m, j$ &
  $3, 3, 1$ &
  $3, 3, 2$ &
  $3, 3, 3$ &
  $3, 4, 1$ \\*[-1mm]
  $Z \left( n, m, j \right)$ &
  $\left[ 486, 28 \right]$ &
  $\left[ 972, 31 \right]$ &
  $\left[ 1944, 37 \right]$ &
  $\left[ 1458, 618 \right]$ \\*[-1mm]
  $Z' \left( n, m, j \right)$ &
  $\left[ 486, 26 \right]$ &
  $\left[ 972, 29 \right]$ &
  $\left[ 1944, 35 \right]$ &
  $\left[ 1458, 615 \right]$ \\ \hline
  $n, m, j$ &
  $6, 2, 1$ &
  $6, 2, 2$ &
  $6, 3, 1$ &
  \\*[-1mm]
  $Z \left( n, m, j \right)$ &
  $\left[ 648, 260 \right]$ &
  $\left[ 1296, 689 \right]$ &
  $\left[ 1944, 833 \right]$ &
  \\*[-1mm]
  $Z' \left( n, m, j \right)$ &
  $D_{18,6}^{(1)}$ &
  $\left[ 1296, 688 \right]$ &
  $\left[ 1944, 832 \right]$ &
  \\ \hline
\end{tabular}
\end{center}
\captionsetup{width=11.5cm}
\caption{The {\tt SmallGroups} identifiers
  of the groups $Z \left( n, m, j \right)$ and $Z' \left( n, m, j \right)$
  with order smaller than 2\,000.
\label{gg}}
\end{table}
The groups $Z \left( n, m, j \right)$ and $Z' \left( n, m, j \right)$
are generated by the same matrices as the groups
$Z' \left( n, m \right)$ and $Z'' \left( n, m \right)$,
respectively,
with the addition of the further generator $- F_{1,j}$,
where $F_{m,j}$ is given in equation~\eqref{xi}.
Notice that there are no groups $Z' \left( n, 2, 1 \right)$ in table~\ref{gg},
because all the matrices generating $Z' \left( n, 2, 1 \right)$,
\textit{viz.}\ $E$,
$L_n$,
$X_2(2)$,
and $- F_{1,1}$ have unit determinant
and therefore $Z^\prime \left( n, 2, 1 \right)$
is a subgroup of $SU(3)$.

The groups $Z \left( n, m, j \right)$ and $Z' \left( n, m, j \right)$
have the same numbers of irreps of each dimension:
$3^{m-1}\, 2^j$ inequivalent singlet irreps,
$3^{m-1}\, 2^{j-1}$
inequivalent doublet irreps,
$\left( n - 1 \right) 3^{m-1}\, 2^j$ triplet irreps,
and $\left( n - 1 \right) \left( n - 2 \right) 3^{m-2}\, 2^{j-2}$
six-plet irreps.

\paragraph{Groups $H \left( n, m, j \right)$:}
When we use generators $E$,
$L_n$,
$X_1(2)$,
and $- F_{m,j}$ with $m>1$,
we obtain groups that we call $H \left( n, m, j \right)$
and list in table~\ref{dugio}.
\begin{table}
\begin{center}
\renewcommand{\arraystretch}{1.2}
\begin{tabular}{|c||c|c|c|c|c|} \hline
  $n, m, j$ &
  $3, 2, 1$ &
  $3, 2, 2$ &
  $3, 2, 3$ &
  $3, 3, 1$ &
  $6, 2, 1$ \\*[-1mm]
  identifier &
  $\left[ 486, 125 \right]$ &
  $\left[ 972, 309 \right]$ &
  $\left[ 1944, 707 \right]$ &
  $\left[ 1458, 1095 \right]$ &
  $\left[ 1944, 2363 \right]$ \\ \hline
\end{tabular}
\end{center}
\captionsetup{width=13.5cm}
\caption{The {\tt SmallGroups} identifiers
  of the groups $H \left( n, m, j \right)$ of order smaller than 2\,000.
\label{dugio}}
\end{table}
The groups $H \left( n,  m, 1 \right)$ have structure $\left[  \left(
  \mathbbm{Z}_{3^{m-1}n} \times \mathbbm{Z}_n 
  \times \mathbbm{Z}_3 \right) \rtimes \mathbbm{Z}_{3}
\right] \rtimes \mathbbm{Z}_2$
and order $3^{m+1}\times2n^2$.
The groups $H \left( n,  m, j \right)$ with $j > 1$
are described in the paragraph of groups $G \left( m, j \right)$ below.

The groups $H \left( n, m, j \right)$ have exactly the same number
of inequivalent irreps of each dimension as the groups
$Z \left( n, m+1, j \right)$ and $Z' \left( n, m+1, j \right)$.

\paragraph{Groups $Y \left( m, j \right)$:}
The groups $Y \left( m, j \right)$,
where $m \ge 2$ and $j \ge 1$,
have structure $\left[ \left( \mathbbm{Z}_{2^j} \times \mathbbm{Z}_{2^j}
  \right) \rtimes \mathbbm{Z}_{3^{m+1}} \right] \rtimes \mathbbm{Z}_3$
and order $3^{m+2}\,4^j$.
There are only three groups $Y \left( m, j \right)$
of order smaller than 2\,000:
\bs
\ba
Y \left( 2, 1 \right) & \cong& \left[ 324, 45 \right], \\
Y \left( 3, 1 \right) & \cong& \left[ 972, 147 \right], \\
Y \left( 2, 2 \right) & \cong& \left[ 1296, 222 \right].
\ea
\es
The groups $Y \left( m, j \right)$
are generated by $L_3$ in equation~\eqref{Pgen2},
$L_{2^j}$ in equation~\eqref{nu}, and
$Z_m$ in equation~\eqref{Pgen3}.

The groups $Y \left( m, j \right)$ only have singlet and triplet irreps:
$3^{m+1}$ inequivalent singlets
and $3^m\, 4^j - 3^{m-1}$ inequivalent triplets.

\paragraph{Groups $G \left( m, j \right)$ and $\left[ 1296,\, 699 \right]$:}
The groups $G \left( m, j \right)$,
where $m \ge 1$ and $j \ge 2$,
have structure $\left\{ \left[ \left(
  \mathbbm{Z}_{3^m} \times \mathbbm{Z}_3 \right)
  \rtimes \mathbbm{Z}_3 \right] \rtimes \mathbbm{Z}_{2^j}
\right\} \rtimes \mathbbm{Z}_3$
and order $3^{m+3}\,2^j$.
The groups $G \left( m, j \right)$ of order smaller than 2\,000
are in table~\ref{G-III}.
\begin{table}
\begin{center}
\renewcommand{\arraystretch}{1.2}
\begin{tabular}{|c||c|c|c|c|c|} \hline
  $ m, j $ &
  $ 1, 2 $ &
  $ 2, 2 $ &
  $ 1, 3 $ &
  $ 2, 3 $ &
  $ 1, 4 $ \\*[-1mm]
  identifier &
  $\left[ 324, 13 \right]$ &
  $\left[ 972, 309 \right]$ &
  $\left[ 648, 19 \right]$ &
  $\left[ 1944, 707 \right]$ &
  $\left[ 1296, 35 \right]$ \\ \hline
\end{tabular}
\end{center}
\captionsetup{width=12.2cm}
\caption{The {\tt SmallGroups} identifiers
  of the groups $G \left( m, j \right)$ with order smaller than 2\,000.
\label{G-III}}
\end{table}
(Notice that the groups $\left[ 972, 309 \right]$ and $\left[ 1944, 707 \right]$
appear in table~\ref{dugio} too.)
The groups $G \left( m, j \right)$ are generated by the matrices $E$,
$- F_{m,j}$,
where $F_{m,j}$ is given in equation~\eqref{xi},
and $\mathrm{diag} \left( 1,\, 1,\, \omega \right)$.
For $m=1$ and $j=2$ one may add a fourth generator $L_2$,
given in equation~\eqref{a4},
to obtain the group $\left[ 1296,\, 699 \right]$,
which has structure $\left\{ \left[ \left(
  \mathbbm{Z}_{6} \times \mathbbm{Z}_6 \right)
  \rtimes \mathbbm{Z}_3 \right] \rtimes \mathbbm{Z}_{4}
\right\} \rtimes \mathbbm{Z}_3$.

The groups $G \left( m, j \right)$ have exactly the same number
of inequivalent irreps of each dimension as the groups
$Z \left( 3, m+1, j \right)$ and $Z' \left( 3, m+1, j \right)$.

\paragraph{Groups $Y(j)$ and $\tilde Y(j)$:}
The groups $Y(j)$ have order $81\times4^j$ and structure
$\left( \mathbbm{Z}_{3\times2^j} \times \mathbbm{Z}_{3\times2^j}
\times \mathbbm{Z}_{3} \right) \rtimes \mathbbm{Z}_3$.
There are three groups $Y(j)$ with order smaller than 2\,000:
\bs
\ba
& & \left[ 81, 7 \right]\ \mbox{with}\ j=0, \\
& & \left[ 324, 60 \right]\ \mbox{with}\ j=1, \\
& & \left[ 1296, 237 \right]\ \mbox{with}\ j=2.
\ea
\es
The group $Y(0)$ coincides with the group
$\Sigma (81)$ or $\Sigma \left( 3\times3^3 \right)$ of ref.~\cite{sigma81}. 
The generators of  $Y(j)$
are the matrix $E$ of equation~\eqref{E} together with the matrix
\be
\label{vugpi}
\mathrm{diag} \left( \xi,\, \xi,\, \xi^2 \right) \quad
\mbox{where} \
\xi = \exp{\left[ 2 i \pi \left/\, \left(3\times2^j \right) \right.\right]}.
\ee
The groups $\tilde Y(j)$ have structure
$\left[ \left(
  \mathbbm{Z}_{3\times2^j} \times \mathbbm{Z}_{3\times2^j} \times \mathbbm{Z}_{3}
\right) \rtimes \mathbbm{Z}_3 \right] \rtimes \mathbbm{Z}_2$
and order $162\times4^j$.
There are two groups $\tilde Y(j)$
with order smaller than 2\,000:
\bs
\ba
& & \left[ 162, 10 \right]\ \mbox{with}\ j=0, \\
& & \left[ 648, 266 \right]\ \mbox{with}\ j=1.
\ea
\es
The generators of  $\tilde{Y}(j)$
are those of $Y(j)$
together with the additional matrix $I'$ in equation~\eqref{Igen2}.

The groups $Y(j)$ have nine inequivalent singlet irreps;
all their remaining irreps are triplets.
The groups $\tilde Y (j)$ have six singlet and three doublet irreps;
their remaining irreps are 12 triplets and one six-plet for $\tilde Y (0)$,
30 triplets and ten six-plets for $\tilde Y (1)$.

\paragraph{Groups $U \left( n, m, j \right)$:}
The groups $U \left( n, m, j \right)$,
where $n$ is a multiple of 3,
$m>1$,
and $1 < j \le m$,
have structure
$\left( \mathbbm{Z}_{3^{m-1}n} \times \mathbbm{Z}_{n} \times \mathbbm{Z}_{3}
\right) \rtimes \mathbbm{Z}_3$ and order $3^{m+1} n^2$. 
We have found the following groups $U \left( n, m, j \right)$
with order smaller than 2\,000:
\bs
\label{buifp}
\ba
& & \left[ 243,\, 55 \right]\ \mbox{with}\ n=3,\ m=2,\ j=2, \\
& & \left[ 729,\, 86 \right]\ \mbox{with}\ n=3,\ m=3,\ j=2, \\
& & \left[ 729,\, 284 \right]\ \mbox{with}\ n=3,\ m=3,\ j=3, \\
& & \left[ 972,\, 550 \right]\ \mbox{with}\ n=6,\ m=2,\ j=2.
\ea
\es
The generators of  $U \left( n, m, j \right)$
are the matrix~$E$ together with
\be
\mathrm{diag} \left( \nu,\, \nu,\, \nu^2 \right),
\quad
\mbox{where}\ \nu = \exp{\left( 2 i \pi \left/\, n \right.\right)},
\ee
and
\be
\label{gihjp1}
\mu\ T_1(m - j + 1),
\quad \mbox{where}\
\mu = \exp{\left[ 2 i \pi \left/\, \left( 3^m \right) \right. \right]}.
\ee
Notice that,
when $j = m$---this happens in three out of the four groups
$U \left( n, m, j \right)$ in~\eqref{buifp}---the
matrix~\eqref{gihjp1} reduces to the matrix
$Y_1(m)$
in equation~\eqref{y1}.

The groups $U \left( n, m, j \right)$ possess
$3^{j+1}$ inequivalent singlet irreps;
all their remaining irreps are triplets.

\paragraph{Groups $L(m)$ and $\left[ 1701, 102 \right]$:}
The groups $L(m)$ have order $3^{m+3}$ and structure
$\left( \mathbbm{Z}_{3^{m+1}} \rtimes \mathbbm{Z}_3 \right)
\rtimes \mathbbm{Z}_3$.
They are generated  by the matrices $X_1(2)$,
$Z_m$,
and $L_3$.
There are the following groups $L(m)$ of order smaller than 2\,000:
\bs
\ba
L(2) &\cong& \left[ 243, 16 \right], \\
L(3) &\cong& \left[ 729, 62 \right].
\ea
\es
Adding the matrix $B_{7,2}$ to the matrices $Z_2$,
$X_1(2)$,
and $L_3$ one generates a group with
structure $\left[ \left( \mathbbm{Z}_{7} \rtimes \mathbbm{Z}_{27} \right)
\rtimes \mathbbm{Z}_3  \right]
\rtimes \mathbbm{Z}_3$ and
{\tt SmallGroups} identifier $\left[1701, 102 \right]$.

The groups $L(m)$ have $3^{m+1}$ singlet irreps;
their other irreps are triplets.

\paragraph{Groups $V(j)$:}
The groups $V(j)$ have order $81\times4^j$ and structure
\be
\left(\mathbbm{Z}_{2^j} \times \mathbbm{Z}_{2^j} \right) \rtimes
\left\{\mathbbm{Z}_{3}\, . \left[\left(\mathbbm{Z}_{3}
  \times \mathbbm{Z}_{3} \right) \rtimes \mathbbm{Z}_3 \right]
= \left(\mathbbm{Z}_{3} \times \mathbbm{Z}_{3} \right)
. \left(\mathbbm{Z}_{3} \times \mathbbm{Z}_{3} \right)\right\}. 
\ee
There are three groups $V(j)$ with order smaller than 2\,000:
\bs
\ba
V(0) &\cong& \left[ 81, 10 \right], \\
V(1) &\cong& \left[ 324, 51 \right], \\
V(2) &\cong& \left[ 1296, 226 \right].
\ea
\es
The generators of $V(j)$ are the matrices $Z_1$,
$X_2 (2)$,
and $L_{2^j}$.

The groups $V(j)$ have nine singlet irreps.
All their other irreps are triplets.

\paragraph{Groups $D(j)$:}
The groups $D(j)$ have structure
$\left( \mathbbm{Z}_{9 \times 2^j} \times \mathbbm{Z}_{9 \times 2^j} \right)
\rtimes \mathbbm{Z}_3$ and order $243 \times 4^j$.
They are generated  by the matrices $E_2$,
$L_{2^j}$, and $T_1(2)$.
There are two groups of order smaller than 2\,000:
\bs
\ba
D(0) &\cong& \left[ 243, 25 \right], \\
D(1) &\cong& \left[ 972, 121 \right].
\ea
\es
Both these groups have nine inequivalent singlets;
their other irreps are triplets.

\paragraph{Groups $J(m)$:}
The groups $J(m)$ have structure
$\mathbbm{Z}_{3^m}\, . \left[
  \left( \mathbbm{Z}_{9} \times \mathbbm{Z}_{3} \right)
\rtimes \mathbbm{Z}_3 \right]$ and order $81\times3^m$.
They are generated  by the matrices $Z_m$ and
$L_9$.
There are two groups of order smaller than 2\,000:
\bs
\ba
J(1) &\cong& \left[ 243, 27 \right], \\
J(2) &\cong& \left[ 729, 80 \right].
\ea
\es
Notice that $J(1)$ coincides with $X(9)$ in table~\ref{Xn}.
The groups $J(m)$ have $3^{m+1}$ singlets;
their other irreps are triplets.

\subsection{The generators of a few more groups} \label{U3sec5}

In this subsection we collect a few more groups together with their generators.

\paragraph{Three groups of order 729:}
Both groups $\left[ 729, 97 \right]$ and $\left[ 729, 98 \right]$
have structure $\left( \mathbbm{Z}_{27} \times \mathbbm{Z}_{9} \right)
\rtimes \mathbbm{Z}_3$.
Group $\left[ 729, 96 \right]$ has the more complicated structure
$\mathbbm{Z}_{3}\, . \left[ \left( \mathbbm{Z}_{9} \times \mathbbm{Z}_{9} \right)
\rtimes \mathbbm{Z}_3 \right]$.
They are generated by the matrix $Z_1$ together with
\bs
\label{Jgen}
\ba
\hat \mu\ \mathrm{diag} \left( \tilde \mu^2,\, \omega,\, \omega \right)
& & \mbox{for}\, \left[ 729, 96 \right], \label{Jgen96}
\\
\hat \mu\ \mathrm{diag} \left( \omega,\, \omega \tilde \mu,\, \omega \tilde \mu
\right)
& & \mbox{for}\, \left[ 729, 97 \right], \label{Jgen97}
\\
\hat \mu\ \mathrm{diag} \left( \omega^2,\,
\omega \tilde \mu,\, \omega \tilde \mu
\right)
& & \mbox{for}\, \left[ 729, 98 \right], \label{Jgen98}
\ea
\es
where
$\omega = \exp{\left( 2 i \pi / 3 \right)}$,
$\tilde \mu = \exp{\left( 2 i \pi / 9 \right)}$,
and $\hat \mu = \exp{\left( 2 i \pi / 27 \right)}$.

Each of these three groups of order 729 possesses
nine singlet and 80 triplet irreps.

\paragraph{A group of order 972:}
The group $\left[ 972, 170 \right]$ is generated by the matrices $L_2$,
$Z_2$,
and $\mathrm{diag} \left( 1,\, 1,\, \omega \right)$.
It has structure
$\left\{ \left[ \left( \mathbbm{Z}_2 \times \mathbbm{Z}_2 \right)
  \rtimes \mathbbm{Z}_{27} \right] \rtimes \mathbbm{Z}_3 \right\}
\rtimes \mathbbm{Z}_3$.

\paragraph{Two groups of order 1458:}
The groups $\left[ 1458, 663 \right]$ and $\left[ 1458, 666 \right]$
have structure $\left[ \left( \mathbbm{Z}_{27} \times \mathbbm{Z}_{9} \right)
  \rtimes \mathbbm{Z}_3 \right] \rtimes \mathbbm{Z}_2$.
They are generated by the matrices $E$,
$L_3$,
$I$,
and~\eqref{Jgen97} for group $\left[ 1458, 663 \right]$,
\eqref{Jgen98} for group $\left[ 1458, 666 \right]$.\footnote{Using
  \eqref{Jgen96} leads to the group
  $\left[ 1458, 659 \right] \cong D_{27,9}^{(1)}$,
  which is subgroup of $SU(3)$.}
Each of these groups has six singlets,
three doublets,
48 triplets,
and 28 six-plets.

\paragraph{Three groups of order 1701:}
The groups with {\tt SmallGroups} identifiers $\left[ 1701, 112 \right]$,
$\left[ 1701, 130 \right]$,
and $\left[ 1701, 131 \right]$ have structure
$\left( \mathbbm{Z}_{7\times3^m} \times \mathbbm{Z}_{3^n} \right)
\rtimes \mathbbm{Z}_3$,
where
$m=3$ and $n=1$ for $\left[ 1701, 112 \right]$
and $m=n=2$
for $\left[ 1701, 130 \right]$ and $\left[ 1701, 131 \right]$.
They are generated by the matrices $B_{7,2}$,
$E$,
$X_2(m)$, and

\bs
\ba
& & T_1(n)\ \mbox{for}\ \left[ 1701, 112 \right]\,
\mbox{and}\, \left[ 1701, 130 \right],
\\
& & T_2(n)\ \mbox{for}\, \left[ 1701, 131 \right].
\ea
\es

\subsection{Other finite subroups of $U(3)$} 
\label{U3sec6}

It is clear from the forms of the generators of
the groups of matrices in subsections~\ref{U3sec2}--\ref{U3sec5}
that those groups are formed by matrices that are all of one of the forms
$R \left( n, a, b, c \right)$,
$V \left( n, a, b, c \right)$,
$W \left( n, a, b, c \right)$ in equations~\eqref{sutyo},
or,
possibly,
also of the forms~\cite{ludl}
\bs
\label{sutyo2}
\ba
S \left( n, a, b, c \right) &\equiv&
\left( \begin{array}{ccc}  \nu^a & 0 & 0 \\ 0 & 0 & \nu^b\\
  0 & \nu^c & 0 \end{array} \right), \\
T \left( n, a, b, c \right) &\equiv&
\left( \begin{array}{ccc}  0 & 0 & \nu^a \\ 0 & \nu^b & 0 \\
  \nu^c & 0 & 0 \end{array} \right), \\
U \left( n, a, b, c \right) &\equiv&
\left( \begin{array}{ccc}  0 & \nu^a & 0 \\ \nu^b & 0 & 0 \\
  0 & 0 & \nu^c \end{array} \right),
\ea
\es
where $\nu = \exp{\left( 2 i \pi / n \right)}$,
for some value of $n$.
It can also be seen that groups RVW only have singlet and triplet irreps,
while groups that include matrices of types~\eqref{sutyo2}
also have doublet and six-plet irreps;
no group in subsections~\ref{U3sec2}--\ref{U3sec5}
has irreps of any other dimension.
There are,
however,
finite subgroups of $U(3)$
that possess no three-dimensional faithful irrep consisting solely
of matrices of the forms~\eqref{sutyo} and~\eqref{sutyo2}.
Those groups have irreps of other dimensions than just one,
two,
three,
and six;
they are analogous to the exceptional subgroups of $SU(3)$.
We present in this section the $U(3)$ subgroups of that type
that have order smaller that 2\,000.

\paragraph{Groups $\Xi \left( m, j \right)$
  and $\hat \Xi \left( m, j \right)$:}
The groups $\Xi\left( m, j \right)$,
where $m \ge 1$ and $j \ge 2$,
have structure
$\left[ \left( \mathbbm{Z}_{3^m} \times \mathbbm{Z}_3 \right)
    \rtimes \mathbbm{Z}_3 \right] \rtimes \mathbbm{Z}_{2^j}$
and order $3^{m+2}\,2^j$.
They are generated by the matrices $E$ in equation~\eqref{E} and $i Q_{m,j}$,
where $Q_{m, j}$ is the matrix in equation \eqref{genQ}.
The groups $\Xi \left( m, j \right)$ of order smaller than 2\,000
are in table~\ref{Xi1}.
Notice that,
since $\det{\left( i Q_{m,j} \right)} = - i \exp{\left[ 2 i \pi
    \left/\, \left( 3^{m-1} 2^j \right) \right. \right]}$,
the group $\Xi \left( 1, 2 \right) \cong \left[ 108, 15 \right]
\cong \Sigma \left( 36\times3 \right)$
is a subgroup of $SU(3)$.
%
\begin{table}
\begin{center}
\renewcommand{\arraystretch}{1.2}
\begin{tabular}{|c||c|c|c|c|c|} \hline
$m, j$ & 1, 2 & 1, 3 & 1, 4 & 1, 5 & 1, 6 \\*[-1mm]
  identifier &
  $\Sigma \left( 36\times3 \right)$ &
  $\left[ 216, 25 \right]$ &
  $\left[ 432, 57 \right]$ &
  $\left[ 864, 194 \right]$ &
  $\left[ 1728, 953 \right]$ \\ \hline
  $m, j$ & 2, 2 & 2, 3 & 2, 4 & 3, 2 & 3, 3 \\*[-1mm]
  identifier &
  $\left[ 324, 111 \right]$ &
  $\left[ 648, 352 \right]$ &
  $\left[ 1296, 1239 \right]$ &
  $\left[ 972, 411 \right]$ &
  $\left[ 1944, 1123 \right]$ \\ \hline  
\end{tabular}
\end{center}
\captionsetup{width=13.6cm}
\caption{The {\tt SmallGroups} identifiers
  of the groups $\Xi\left( m, j \right)$ 
  with order smaller than 2\,000.
\label{Xi1}}
\end{table}

The groups $\hat \Xi \left( m, j \right)$ have structure
$\left\{ \left[ \left( \mathbbm{Z}_{3^m} \times \mathbbm{Z}_3 \right)
  \rtimes \mathbbm{Z}_3 \right] \rtimes \mathbbm{Z}_{2^j}
\right\} \rtimes \mathbbm{Z}_2$
and order $3^{m+2}\,2^{j+1}$.
The generators of $\hat \Xi \left( m, j \right)$
are the same as the generators of $\Xi \left( m, j \right)$
together with the additional matrix $I'$ in equation~\eqref{Igen2}.
The groups $\hat \Xi \left( m, j \right)$ of order smaller than 2\,000
are in table~\ref{Xi1hat}.
\begin{table}
\begin{center}
\renewcommand{\arraystretch}{1.2}
\begin{tabular}{|c||c|c|c|c|} \hline
  $m, j$  &
  1, 3 &
  1, 4 &
  1, 5 &
  2, 3 \\*[-1mm]
identifier &
$\left[ 432, 273 \right]$ &
$\left[ 864, 737 \right]$ &
$\left[ 1728, 2929 \right]$ &
$\left[ 1296, 2203 \right]$  \\ \hline
\end{tabular}
\end{center}
\captionsetup{width=11.1cm}
\caption{The {\tt SmallGroups} identifiers
  of the groups $\hat \Xi \left( m, j \right)$ with order smaller than 2\,000.
\label{Xi1hat}}
\end{table}
Notice that
the groups $\hat \Xi \left( m, 2 \right)$ have structure
$\left\{ \left[ \left( \mathbbm{Z}_{3^m} \times \mathbbm{Z}_3 \right)
\rtimes \mathbbm{Z}_3 \right] \rtimes \mathbbm{Z}_4 \right\}
\times \mathbbm{Z}_2$;
these groups may be written as \emph{direct}\/ products of $\mathbbm{Z}_2$
and smaller groups,
hence they are not included in table~\ref{Xi1hat}.

The groups $\Xi \left( m, j \right)$ possess
$3^{m-1}2^j$ inequivalent singlet irreps,
$3^{m-1}2^{j+1}$ inequivalent triplet irreps,
and $3^{m-1}2^{j-1}$ inequivalent quadruplet irreps.
The groups $\hat \Xi \left( m, j \right)$ have twice as many irreps
of each dimension as the groups $\Xi \left( m, j \right)$.

\paragraph{Groups $\Pi \left( m, j \right)$:}
The groups $\Pi \left( m, j \right)$,
where $m \ge 1$ and $j \ge 2$,
have structure  $\left\{ \left[ \left(
  \mathbbm{Z}_{3^m} \times \mathbbm{Z}_3 \right) \rtimes \mathbbm{Z}_3 \right]
\rtimes \mathbbm{Z}_4 \right\} \rtimes \mathbbm{Z}_{2^j}$
and order $3^{m+2}\,2^{j+2}$.
They are generated by the matrices $E$,
$K$,
and $Q_{m,j}$.
The groups $\Pi \left( m, j \right)$ of order smaller than 2\,000
are in table~\ref{Xi2}.
\begin{table}
\begin{center}
\renewcommand{\arraystretch}{1.2}
\begin{tabular}{|c||c|c|c|c|} \hline
  $m, j$  &
  1, 2 &
  1, 3 &
  1, 4 &
  2, 2 \\*[-1mm]
  identifier &
  $\left[ 432, 239 \right]$ &
  $\left[ 864, 675 \right]$ &
  $\left[ 1728, 2785 \right]$ &
  $\left[ 1296, 1995 \right]$  \\ \hline
\end{tabular}
\end{center}
\captionsetup{width=11.1cm}
\caption{The {\tt SmallGroups} identifiers
  of the groups $\Pi \left( m, j \right)$ with order smaller than 2\,000.
\label{Xi2}}
\end{table}

The groups $\Pi \left( m, j \right)$ have $3^{m-1}2^{j+1}$ singlet irreps,
$3^{m-1}2^{j-1}$ doublet irreps,
$3^{m-1}2^{j+2}$ triplet irreps,
$3^{m-1}2^j$ six-plet irreps,
and $3^{m-1}2^{j-1}$ eight-plet irreps.

\paragraph{Groups $\Theta (m)$:}
These groups have structure
$\left[ \left( \mathbbm{Z}_{3^m} \times \mathbbm{Z}_3 \right)
  \rtimes \mathbbm{Z}_3 \right] \rtimes Q_8$,
where $Q_8$ is the quaternion group.
Since $Q_8$ has eight elements,
the order of $\Theta (m)$ is $72\times3^m$.
The generators of $\Theta (m)$ are $E$,
$K$,
and either $Q_{m,1}$ or $Q_{m,0}$.
Since $\det{Q_{m,0}} = \exp{\left[ 2 i \pi
    \left/\, \left( 3^{m-1} \right) \right. \right]}$,
the group
$\Theta (1)$
is a subgroup of $SU(3)$.
Notice that the groups $\Theta(m)$ have the same generators
as hypothetical groups $\Pi \left( m, 1 \right)$ would have had;
but they have a slightly different structure.
There are three groups $\Theta (m)$ of order smaller than 2\,000:
\bs
\ba
\Theta (1) &\cong& \Sigma \left( 72 \times 3 \right), \\
\Theta (2) &\cong& \left[ 648, 551 \right], \\
\Theta (3) &\cong& \left[ 1944, 2333 \right].
\ea
\es
The groups $\Theta (m)$ have
as many inequivalent irreps of each dimension
as groups $\Pi \left( m, 1 \right)$.

\paragraph{Groups $\Upsilon (m)$ and $\Upsilon'(m)$:}
These groups have structure
\be
\left\{ \left[ \left( \mathbbm{Z}_3 \times \mathbbm{Z}_3 \right)
  \rtimes \mathbbm{Z}_3 \right]
\rtimes Q_8 \right\} .\, \mathbbm{Z}_{3^{m-1}}
= \mathbbm{Z}_{3^{m-1}} \, .\left\{ \left[ \left(
    \mathbbm{Z}_3 \times \mathbbm{Z}_3 \right)
    \rtimes Q_8 \right] \rtimes \mathbbm{Z}_3 \right\} 
  \ee
and order $72\times3^m$.
The generators are $E$,
$Q_{0,0}$,
and $X_1(m)$ for $\Upsilon(m)$ or $X_2(m)$ for $\Upsilon'(m)$.
There are the following groups of order smaller than 2\,000:
\bs
\ba
\Upsilon(2) &\cong& \left[ 648, 531 \right], \\
\Upsilon (3) &\cong& \left[ 1944, 2293 \right], \\
\Upsilon'(2) &\cong& \Sigma \left( 216\times3 \right), \\
\Upsilon'(3) &\cong& \left[ 1944, 2294 \right].
\ea
\es
Notice that all three generators of $\Upsilon'(2)$
have unit determinant and therefore $\Upsilon'(2)$
is a subgroup of $SU(3)$.

The groups $\Upsilon (m)$ and $\Upsilon' (m)$ have
$3^{m-1}$ singlets,
$3^{m-1}$ doublets,
$7\times3^{m-2}$ triplets,
$2\times3^{m-1}$ six-plets,
$3^{m-1}$ eight-plets,
and $3^{m-2}\times2$ nine-plets.

\paragraph{Groups $\Omega (m)$:}
These groups have structure
$\left\{ \left[ \left( \mathbbm{Z}_{3^m} \times \mathbbm{Z}_3 \right)
  \rtimes \mathbbm{Z}_3 \right] \rtimes Q_8 \right\} \rtimes \mathbbm{Z}_3$
and order $72\times3^{m+1}$.
They are generated by the matrices $Q_{m,0}$ and $Z_1$. 
There are the following groups of order smaller than 2\,000:
\bs
\ba
\Omega (1) &\cong& \left[ 648, 533 \right], \\
\Omega (2) &\cong& \left[ 1944, 3448 \right].
\ea
\es
The groups $\Omega (m)$ have exactly as many inequivalent irreps
of each dimension as the groups $\Upsilon (m+1)$ and $\Upsilon' (m+1)$.

\section{Conclusion}

In this paper we have used the {\tt SmallGroups} library
to search for all the finite subgroups of $U(3)$
of order less than 2\,000
that have a faithful three-dimensional irreducible representation
and that cannot be written as the direct product
of some smaller group and a cyclic group.
We have found that there are three types of finite subgroups of $U(3)$:
\begin{itemize}
\item Groups that have a three-dimensional representation
  consisting solely of matrices of the forms~\eqref{sutyo}
  for some value of $n$.
  Those groups only have singlet and triplet irreducible representations.
\item Groups that have a three-dimensional representation
  consisting solely of matrices of the forms~\eqref{sutyo} and~\eqref{sutyo2}
  for some value of $n$.
  Those groups only have singlet,
  doublet,
  triplet,
  and six-plet irreducible representations.
\item Groups that do \emph{not}\/ have a three-dimensional representation
  consisting solely of matrices of the forms~\eqref{sutyo} and~\eqref{sutyo2}.
  Those groups have irreducible representations of other dimensions,
  like for instance 4-plets,
  8-plets,
  or 9-plets.
  Their generators include matrices $Q_{m,j}$ and possibly $K$
  in equations~\eqref{genK},
  \eqref{genQ}.
  These groups include as special cases the exceptional $SU(3)$ subgroups
  $\Sigma \left( 36\times3 \right)$,
  $\Sigma \left( 72\times3 \right)$,
  and $\Sigma \left( 216\times3 \right)$.\footnote{It seems likely to us
that the $SU(3)$ subgroup $\Sigma \left( 360\times3 \right)$ is also
a special case of a series of $U(3)$ subgroups; the other groups of
that series,
though,
surely have order larger than 2\,000.}
\end{itemize}
We were able to group most finite subgroups of $U(3)$ in many series
depending on one,
two,
or sometimes three integers;
the groups in each series have related generators
and related numbers of irreps of each dimension.
Unfortunately,
many of these series have very few groups
and we do not know whether and how they extend to groups of order
higher than 2\,000.
It is possible
(and it would be desirable)
that some of these series may be further unified
among themselves.

\paragraph{Acknowledgements:}

D.J.~thanks the Lithuanian Academy of Sciences
for support through the project~DaFi2016.
The work of L.L.\ is supported by the Portuguese
\textit{Fun\-da\-\c c\~ao pa\-ra a Ci\-\^en\-cia e a Te\-cno\-lo\-gia}\/
through the projects CERN/FIS-NUC/0010/2015
and UID/FIS/00777/2013,
which are partially funded by POCTI (FEDER),
COMPETE,
QREN,
and the European Union.

\appendix

\section{Full tables}
  
  In this appendix we present tables of all the groups of order
  smaller than 2\,000 that have a faithful three-dimensional irrep
  and cannot be written as the direct product of some smaller group
  and a cyclic group.
  The groups are ordered according to increasing values of firstly $o$
  and then $j$ in their {\tt SmallGroups} identifier $\left[ o, j \right]$.
  Tables~\ref{tab:SU3_par1}--\ref{tab:SU3_par7} include the groups
  that are subgroups of $SU(3)$;
  tables~\ref{tab:U3_par1}--\ref{tab:U3_par11} include the groups
  that are not subgroups of $SU(3)$.

\begin{table}
\begin{centering}
\begin{minipage}{0.46\textwidth}
\renewcommand{\arraystretch}{1.15}
\begin{tabular}{|l|l|}\hline
Identifier & Classification \\ \hline \hline 
$\left[ 12, 3\right]$ & $ \Delta \left( 3\times2^2 \right)$ \\
$\left[ 21, 1\right]$ & $ C_{7,1}^{(2)}$ \\
$\left[ 24, 12 \right]$ & $ \Delta \left( 6\times2^2 \right)$ \\
$\left[ 27, 3\right]$ & $ \Delta \left( 3\times3^2 \right)$ \\
$\left[ 39, 1\right]$ & $ C_{13,1}^{(3)}$ \\
$\left[ 48, 3\right]$ & $ \Delta \left( 3\times4^2 \right)$ \\
$\left[ 54, 8\right]$ & $ \Delta \left( 6\times3^2 \right)$ \\
$\left[ 57, 1\right]$ & $ C_{19,1}^{(7)}$ \\
$\left[ 60, 5\right]$ & $ \Sigma \left( 60 \right)$ \\
$\left[ 75, 2\right]$ & $ \Delta \left( 3\times5^2 \right)$ \\
$\left[ 81, 9\right]$ & $ C_{9,3}^{(1)}$ \\
$\left[ 84, 11 \right]$ & $ C_{14,2}^{(2)} $ \\
$\left[ 93, 1\right]$ & $ C_{31,1}^{(5)} $ \\
$\left[ 96, 64 \right]$ & $ \Delta \left( 6\times4^2 \right)$ \\
$\left[ 108, 15\right]$ & $ \Sigma \left( 36\times3 \right)$ \\
$\left[ 108, 22\right]$ & $ \Delta \left( 3\times6^2 \right)$ \\
$\left[ 111, 1 \right]$ & $ C_{37,1}^{(10)}$ \\
$\left[ 129, 1 \right]$ & $ C_{43,1}^{(6)} $ \\
$\left[ 147, 1 \right]$ & $ C_{49,1}^{(18)}$ \\
$\left[ 147, 5 \right]$ & $ \Delta \left( 3\times7^2 \right)$ \\
$\left[ 150, 5 \right]$ & $ \Delta \left( 6\times5^2 \right)$ \\
$\left[ 156, 14\right]$ & $ C_{26,2}^{(3)} $ \\
$\left[ 162, 14\right]$ & $ D_{9,3}^{(1)} $ \\
$\left[ 168, 42\right]$ & $ \Sigma \left( 168 \right)$ \\
$\left[ 183, 1 \right]$ & $ C_{61,1}^{(13)}$ \\
$\left[ 189, 8 \right]$ & $ C_{21,3}^{(2)} $ \\
$\left[ 192, 3 \right]$ & $ \Delta \left( 3\times8^2 \right) $ \\
$\left[ 201, 1 \right]$ & $ C_{67,1}^{(29)}$ \\
\hline
\end{tabular}
\captionsetup{width=6.3cm}
\caption{The finite subgroups of $SU(3)$.
Part~1: groups with order through 201.}
\label{tab:SU3_par1} 
\end{minipage}
\hfill
\begin{minipage}{0.46\textwidth}
\renewcommand{\arraystretch}{1.138}
\begin{tabular}{|l|l|}\hline
Identifier & Classification \\ \hline \hline 
$\left[ 216, 88\right]$ & $ \Sigma \left( 72\times3 \right)$ \\
$\left[ 216, 95\right]$ & $ \Delta \left( 6\times6^2 \right)$ \\
$\left[ 219, 1 \right]$ & $ C_{73,1}^{(8)} $ \\
$\left[ 228, 11\right]$ & $ C_{38,2}^{(7)} $ \\
$\left[ 237, 1 \right]$ & $ C_{79,1}^{(23)}$ \\
$\left[ 243, 26\right]$ & $ \Delta \left( 3\times9^2 \right)$ \\
$\left[ 273, 3 \right]$ & $ C_{91,1}^{(16)}$ \\
$\left[ 273, 4 \right]$ & $ C_{91,1}^{(9)} $ \\
$\left[ 291, 1 \right]$ & $ C_{97,1}^{(35)}$ \\
$\left[ 294, 7 \right]$ & $ \Delta \left( 6\times7^2 \right)$ \\
$\left[ 300, 43\right]$ & $ \Delta \left( 3\times10^2 \right)$ \\
$\left[ 309, 1 \right]$ & $ C_{103,1}^{(46)} $ \\
$\left[ 324, 50\right]$ & $ C_{18,6}^{(1)} $ \\
$\left[ 327, 1 \right]$ & $ C_{109,1}^{(45)} $ \\
$\left[ 336, 57\right]$ & $ C_{28,4}^{(2)} $ \\
$\left[ 351, 8 \right]$ & $ C_{39,3}^{(3)} $ \\
$\left[ 363, 2 \right]$ & $ \Delta \left( 3\times11^2 \right)$ \\
$\left[ 372, 11\right]$ & $ C_{62,2}^{(5)} $ \\
$\left[ 381, 1 \right]$ & $ C_{127,1}^{(19)} $ \\
$\left[ 384, 568 \right]$ & $ \Delta \left( 6\times8^2 \right)$ \\
$\left[ 399, 3 \right]$ & $ C_{133,1}^{(11)} $ \\
$\left[ 399, 4 \right]$ & $ C_{133,1}^{(30)} $ \\
$\left[ 417, 1 \right]$ & $ C_{139,1}^{(42)} $ \\
$\left[ 432, 103 \right]$ & $ \Delta \left( 3\times12^2 \right)$ \\
$\left[ 444, 14\right]$ & $ C_{74,2}^{(10)}$ \\
$\left[ 453, 1 \right]$ & $ C_{151,1}^{(32)} $ \\
$\left[ 471, 1 \right]$ & $ C_{157,1}^{(12)} $ \\
$\left[ 486, 61\right]$ & $ \Delta \left( 6\times9^2 \right)$ \\
\hline
\end{tabular}
\captionsetup{width=6.3cm}
\caption{The finite subgroups of $SU(3)$.
Part~2: groups with $216 \le \mbox{order} \le 486$.}
\label{tab:SU3_par2} 
\end{minipage}
\end{centering}
\end{table}

\begin{table}
\begin{centering}
\begin{minipage}{0.46\textwidth}
\renewcommand{\arraystretch}{1.13}
\begin{tabular}{|l|l|}\hline
Identifier & Classification \\ \hline \hline 
$\left[ 489, 1 \right]$ & $ C_{163,1}^{(58)} $ \\
$\left[ 507, 1 \right]$ & $ C_{169,1}^{(22)} $ \\
$\left[ 507, 5 \right]$ & $ \Delta \left( 3\times13^2 \right)$ \\
$\left[ 513, 9 \right]$ & $ C_{57,3}^{(7)} $ \\
$\left[ 516, 11\right]$ & $ C_{86,2}^{(6)} $ \\
$\left[ 525, 5 \right]$ & $ C_{35,5}^{(2)} $ \\
$\left[ 543, 1 \right]$ & $ C_{181,1}^{(48)} $ \\
$\left[ 567, 13\right]$ & $ C_{63,3}^{(4)} $ \\
$\left[ 579, 1 \right]$ & $ C_{193,1}^{(84)} $ \\
$\left[ 588, 11\right]$ & $ C_{98,2}^{(18)}$ \\
$\left[ 588, 60\right]$ & $ \Delta \left( 3\times14^2 \right)$ \\
$\left[ 597, 1 \right]$ & $ C_{199,1}^{(92)} $ \\
$\left[ 600, 179 \right]$ & $ \Delta \left( 6\times10^2 \right)$ \\
$\left[ 624, 60\right]$ & $ C_{52,4}^{(3)} $ \\
$\left[ 633, 1 \right]$ & $ C_{211,1}^{(14)} $ \\
$\left[ 648, 259 \right]$ & $ D_{18,6}^{(1)}$ \\
$\left[ 648, 532 \right]$ & $ \Sigma \left( 216\times3 \right) $ \\
$\left[ 651, 3 \right]$ & $ C_{217,1}^{(25)} $ \\
$\left[ 651, 4 \right]$ & $ C_{217,1}^{(67)} $ \\
$\left[ 669, 1 \right]$ & $ C_{223,1}^{(39)} $ \\
$\left[ 675, 12\right]$ & $ \Delta \left( 3\times15^2 \right)$ \\
$\left[ 687, 1 \right]$ & $ C_{229,1}^{(94)} $ \\
$\left[ 723, 1 \right]$ & $ C_{241,1}^{(15)} $ \\ 
$\left[ 726, 5 \right]$ & $ \Delta \left( 6\times11^2 \right)$ \\
$\left[ 729, 95\right]$ & $ C_{27,9}^{(1)} $ \\
$\left[ 732, 14\right]$ & $ C_{122,2}^{(13)} $ \\
$\left[ 741, 3 \right]$ & $ C_{247,1}^{(87)} $ \\
$\left[ 741, 4 \right]$ & $ C_{247,1}^{(68)} $ \\
$\left[ 756, 117 \right]$ & $ C_{42,6}^{(2)} $ \\
\hline
\end{tabular}
\captionsetup{width=6.3cm}
\caption{The finite subgroups of $SU(3)$.
Part~3: groups with $489 \le \mbox{order} \le 756$.}
\label{tab:SU3_par3} 
\end{minipage}
\hfill
\begin{minipage}{0.46\textwidth}
\renewcommand{\arraystretch}{1.2}
\begin{tabular}{|l|l|}\hline
Identifier & Classification \\ \hline \hline 
$\left[ 768, 1083477 \right]$ & $ \Delta \left( 3\times16^2 \right)$ \\
$\left[ 777, 3 \right]$ & $ C_{259,1}^{(121)}$ \\
$\left[ 777, 4 \right]$ & $ C_{259,1}^{(100)}$ \\
$\left[ 804, 11\right]$ & $ C_{134,2}^{(29)} $ \\
$\left[ 813, 1 \right]$ & $ C_{271,1}^{(28)} $ \\
$\left[ 831, 1 \right]$ & $ C_{277,1}^{(116)}$ \\
$\left[ 837, 8 \right]$ & $ C_{93,3}^{(5)} $ \\
$\left[ 849, 1 \right]$ & $ C_{283,1}^{(44)} $ \\
$\left[ 864, 701 \right]$ & $ \Delta \left( 6\times12^2 \right)$ \\
$\left[ 867, 2 \right]$ & $ \Delta \left( 3\times17^2 \right)$ \\
$\left[ 876, 14\right]$ & $ C_{146,2}^{(8)}$ \\
$\left[ 903, 5 \right]$ & $ C_{301,1}^{(135)}$ \\
$\left[ 903, 6 \right]$ & $ C_{301,1}^{(79)} $ \\
$\left[ 912, 57\right]$ & $ C_{76,4}^{(7)} $ \\
$\left[ 921, 1 \right]$ & $ C_{307,1}^{(17)} $ \\
$\left[ 939, 1 \right]$ & $ C_{313,1}^{(98)} $ \\
$\left[ 948, 11\right]$ & $ C_{158,2}^{(23)} $ \\
$\left[ 972, 122 \right]$ & $ \Delta \left( 3\times18^2 \right)$ \\
$\left[ 975, 5 \right]$ & $ C_{65,5}^{(3)} $ \\
$\left[ 993, 1 \right]$ & $ C_{331,1}^{(31)} $ \\
$\left[ 999, 9 \right]$ & $ C_{111,3}^{(10)} $ \\
$\left[ 1011, 1\right]$ & $ C_{337,1}^{(128)}$ \\
$\left[ 1014, 7\right]$ & $ \Delta \left( 6\times13^2 \right)$ \\
$\left[ 1029, 6\right]$ & $ C_{343,1}^{(18)} $ \\
$\left[ 1029, 9\right]$ & $ C_{49,7}^{(2)} $ \\
$\left[ 1047, 1\right]$ & $ C_{349,1}^{(122)}$ \\
$\left[ 1053, 35 \right]$ & $ C_{117,3}^{(16)} $ \\
$\left[ 1080, 260\right]$ & $ \Sigma \left( 360\times3 \right) $ \\
\hline
\end{tabular}
\captionsetup{width=6.3cm}
\caption{The finite subgroups of $SU(3)$.
Part~4: groups with $768 \le \mbox{order} \le 1080$.}
\label{tab:SU3_par4} 
\end{minipage}
\end{centering}
\end{table}

\begin{table}
\begin{centering}
\begin{minipage}{0.46\textwidth}
\renewcommand{\arraystretch}{1.15}
\begin{tabular}{|l|l|}\hline
Identifier & Classification \\ \hline \hline 
$\left[ 1083, 1\right]$ & $ C_{361,1}^{(68)} $ \\
$\left[ 1083, 5\right]$ & $ \Delta \left( 3\times19^2 \right)$ \\
$\left[ 1092, 68 \right]$ & $ C_{182,2}^{(9)}$ \\
$\left[ 1092, 69 \right]$ & $ C_{182,2}^{(16)} $ \\
$\left[ 1101, 1\right]$ & $ C_{367,1}^{(83)} $ \\
$\left[ 1119, 1\right]$ & $ C_{373,1}^{(88)} $ \\ 
$\left[ 1137, 1\right]$ & $ C_{379,1}^{(51)} $ \\
$\left[ 1161, 9\right]$ & $ C_{129,3}^{(6)}$ \\
$\left[ 1164, 14 \right]$ & $ C_{194,2}^{(35)} $ \\
$\left[ 1176, 243\right]$ & $ \Delta \left( 6\times14^2 \right)$ \\
$\left[ 1191, 1\right]$ & $ C_{397,1}^{(34)} $ \\
$\left[ 1200, 384\right]$ & $ \Delta \left( 3\times20^2 \right)$ \\
$\left[ 1209, 3\right]$ & $ C_{403,1}^{(87)} $ \\
$\left[ 1209, 4\right]$ & $ C_{403,1}^{(191)}$ \\
$\left[ 1227, 1\right]$ & $ C_{409,1}^{(53)} $ \\
$\left[ 1236, 11 \right]$ & $ C_{206,2}^{(46)} $ \\
$\left[ 1263, 1\right]$ & $ C_{421,1}^{(20)} $ \\
$\left[ 1281, 3\right]$ & $ C_{427,1}^{(135)}$ \\
$\left[ 1281, 4\right]$ & $ C_{427,1}^{(74)} $ \\
$\left[ 1296, 228\right]$ & $ C_{36,12}^{(1)}$ \\
$\left[ 1299, 1\right]$ & $ C_{433,1}^{(198)}$ \\
$\left[ 1308, 14 \right]$ & $ C_{218,2}^{(45)} $ \\
$\left[ 1317, 1\right]$ & $ C_{439,1}^{(171)}$ \\
$\left[ 1323, 8\right]$ & $ C_{147,3}^{(18)} $ \\
$\left[ 1323, 43 \right]$ & $ \Delta \left( 3\times21^2 \right)$ \\
$\left[ 1344, 393\right]$ & $ C_{56,8}^{(2)} $ \\
$\left[ 1350, 46 \right]$ & $ \Delta \left( 6\times15^2 \right)$ \\
$\left[ 1371, 1\right]$ & $ C_{457,1}^{(133)}$ \\
$\left[ 1389, 1\right]$ & $ C_{463,1}^{(21)} $ \\
\hline
\end{tabular}
\captionsetup{width=6.3cm}
\caption{The finite subgroups of $SU(3)$.
Part~5: groups with $1083 \le \mbox{order} \le 1389$.}
\label{tab:SU3_par5} 
\end{minipage}
\hfill
\begin{minipage}{0.46\textwidth}
\renewcommand{\arraystretch}{1.14}
\begin{tabular}{|l|l|}\hline
Identifer & Classification \\ \hline \hline 
$\left[ 1404, 141\right]$ & $ C_{78,6}^{(3)} $ \\
$\left[ 1407, 3\right]$ & $ C_{469,1}^{(163)}$ \\
$\left[ 1407, 4\right]$ & $ C_{469,1}^{(37)} $ \\
$\left[ 1425, 5\right]$ & $ C_{95,5}^{(7)} $ \\
$\left[ 1443, 3\right]$ & $ C_{481,1}^{(100)}$ \\
$\left[ 1443, 4\right]$ & $ C_{481,1}^{(211)}$ \\
$\left[ 1452, 34 \right]$ & $ \Delta \left( 3\times22^2 \right) $ \\
$\left[ 1458, 659\right]$ & $ D_{27,9}^{(1)}$ \\
$\left[ 1461, 1\right]$ & $ C_{487,1}^{(232)}$ \\
$\left[ 1488, 57 \right]$ & $ C_{124,4}^{(5)}$ \\
$\left[ 1497, 1\right]$ & $ C_{499,1}^{(139)}$ \\
$\left[ 1524, 11 \right]$ & $ C_{254,2}^{(19)} $ \\
$\left[ 1533, 3\right]$ & $ C_{511,1}^{(137)}$ \\
$\left[ 1533, 4\right]$ & $ C_{511,1}^{(81)} $ \\
$\left[ 1536, 408544632\right]$ & $ \Delta \left( 6\times16^2 \right)$ \\
$\left[ 1539, 35 \right]$ & $ C_{171,3}^{(7)}$ \\
$\left[ 1569, 1\right]$ & $ C_{523,1}^{(60)} $ \\
$\left[ 1587, 2\right]$ & $ \Delta \left( 3\times23^2 \right)$ \\
$\left[ 1596, 55 \right]$ & $ C_{266,2}^{(11)} $ \\
$\left[ 1596, 56 \right]$ & $ C_{266,2}^{(30)} $ \\
$\left[ 1623, 1\right]$ & $ C_{541,1}^{(129)}$ \\
$\left[ 1641, 1\right]$ & $ C_{547,1}^{(40)} $ \\
$\left[ 1647, 9\right]$ & $ C_{183,3}^{(13)} $ \\
$\left[ 1659, 3\right]$ & $ C_{553,1}^{(23)} $ \\
$\left[ 1659, 4\right]$ & $ C_{553,1}^{(102)}$ \\
$\left[ 1668, 11 \right]$ & $ C_{278,2}^{(42)} $ \\
$\left[ 1677, 3\right]$ & $ C_{559,1}^{(165)}$ \\
$\left[ 1677, 4\right]$ & $ C_{559,1}^{(178)}$ \\
$\left[ 1701, 135\right]$ & $ C_{63,9}^{(2)} $ \\
\hline
\end{tabular}
\captionsetup{width=6.3cm}
\caption{The finite subgroups of $SU(3)$.
Part~6: groups with $1404 \le \mbox{order} \le 1701$.}
\label{tab:SU3_par6} 
\end{minipage}
\end{centering}
\end{table}

\begin{table}
\begin{centering}
\begin{minipage}{0.46\textwidth}
\renewcommand{\arraystretch}{1.2}
\begin{tabular}{|l|l|}\hline
Identifier & Classification \\ \hline \hline 
$\left[ 1713, 1\right]$ & $ C_{571,1}^{(109)}$ \\
$\left[ 1728, 1291 \right]$ & $ \Delta \left( 3\times24^2 \right)$ \\
$\left[ 1731, 1\right]$ & $ C_{577,1}^{(213)}$ \\
$\left[ 1734, 5\right]$ & $ \Delta \left( 6\times17^2 \right)$ \\
$\left[ 1767, 3\right]$ & $ C_{589,1}^{(87)} $ \\
$\left[ 1767, 4\right]$ & $ C_{589,1}^{(273)}$ \\
$\left[ 1776, 60 \right]$ & $ C_{148,4}^{(10)} $ \\
$\left[ 1803, 1\right]$ & $ C_{601,1}^{(24)} $ \\
$\left[ 1809, 9\right]$ & $ C_{201,3}^{(29)} $ \\
$\left[ 1812, 11 \right]$ & $ C_{302,2}^{(32)} $ \\
$\left[ 1821, 1\right]$ & $ C_{607,1}^{(210)}$ \\
$\left[ 1839, 1\right]$ & $ C_{613,1}^{(65)} $ \\
$\left[ 1857, 1\right]$ & $ C_{619,1}^{(252)}$ \\
$\left[ 1875, 16 \right]$ & $ \Delta \left( 3\times25^2 \right)$ \\
$\left[ 1884, 14 \right]$ & $ C_{314,2}^{(12)} $ \\
$\left[ 1893, 1\right]$ & $ C_{631,1}^{(43)} $ \\
$\left[ 1911, 3\right]$ & $ C_{637,1}^{(165)}$ \\
$\left[ 1911, 4\right]$ & $ C_{637,1}^{(263)}$ \\
$\left[ 1911, 14 \right]$ & $ C_{91,7}^{(3)} $ \\
$\left[ 1929, 1\right]$ & $ C_{643,1}^{(177)}$ \\
$\left[ 1944, 849\right]$ & $ \Delta \left( 6\times18^2 \right)$ \\
$\left[ 1956, 11 \right]$ & $ C_{326,2}^{(58)} $ \\
$\left[ 1971, 9\right]$ & $ C_{219,3}^{(8)}$ \\
$\left[ 1983, 1\right]$ & $ C_{661,1}^{(296)}$ \\
\hline
\end{tabular}
\captionsetup{width=6.3cm}
\caption{The finite subgroups of $SU(3)$.
Part~7: groups with $1713 \le \mbox{order} < 2000$.}
\label{tab:SU3_par7} 
\end{minipage}
\end{centering}
\end{table}

\newpage

\begin{table}
\begin{centering}
\begin{minipage}{0.46\textwidth}
\renewcommand{\arraystretch}{1.1}
\begin{tabular}{|l|l|}\hline
Identifier & Classification \\\hline\hline 
$\left[ 27, 4\right]$ & $ X(3),\; W(1,2) $ \\
$\left[ 36, 3\right]$ & $ \Delta \left( 3\times2^2,2 \right)$ \\
$\left[ 48, 30 \right]$ & $ S_4(2) $ \\
$\left[ 63, 1\right]$ & $ T_{7}^{(2)}(2)$ \\
$\left[ 81, 6\right]$ & $ W \left( 1, 3 \right) $ \\
$\left[ 81, 7\right]$ & $ Y(0),\; \Sigma \left( 3\times3^3 \right)$ \\
$\left[ 81, 8\right]$ & $ Z' \left( 3, 2 \right) $ \\
$\left[ 81, 10 \right]$ & $ V \left( 0 \right) $ \\
$\left[ 81, 14 \right]$ & $ Z \left( 3, 2 \right)$ \\
$\left[ 96, 65 \right]$ & $ S_4(3) $ \\
$\left[ 108, 3 \right]$ & $ \Delta \left( 3\times2^2,3 \right)$ \\
$\left[ 108, 11\right]$ & $ \Delta \left( 6\times3^2,2 \right)$ \\
$\left[ 108, 19\right]$ & $ W \left( 2, 2 \right) $ \\
$\left[ 108, 21\right]$ & $ X(6)$ \\
$\left[ 117, 1 \right]$ & $ T_{13}^{(3)}(2)$ \\
$\left[ 144, 3 \right]$ & $ \Delta \left( 3\times4^2,2 \right)$ \\
$\left[ 162, 10\right]$ & $ \tilde Y(0) $ \\
$\left[ 162, 12\right]$ & $ Z \left( 3, 2, 1 \right)$ \\
$\left[ 162, 44\right]$ & $ \Delta ' \left( 6\times3^2,2,1 \right)$ \\
$\left[ 171, 1 \right]$ & $ T_{19}^{(7)}(2)$ \\
$\left[ 189, 1 \right]$ & $ T_{7}^{(2)}(3)$ \\
$\left[ 189, 4 \right]$ & $ Q_{7}^{(2)}(2)$ \\
$\left[ 189, 5 \right]$ & $ Q_{7}^{(2)\prime}(2)$ \\
$\left[ 189, 7 \right]$ & $ P_{7}^{(2)}(2)$ \\
$\left[ 192, 182\right]$ & $ \Delta \left( 6\times4^2,2 \right)$ \\
$\left[ 192, 186\right]$ & $ S_4(4) $ \\
$\left[ 216, 17\right]$ & $ \Delta \left( 6\times3^2,3 \right)$ \\
$\left[ 216, 25\right]$ & $ \Xi \left( 1,3 \right) $ \\
$\left[ 225, 3 \right]$ & $ \Delta \left( 3\times5^2,2 \right)$ \\
\hline
\end{tabular}
\caption{The finite subgroups of $U(3)$.
Part~1: groups with $\mbox{order} \le 225$.}
\label{tab:U3_par1} 
\end{minipage}
\hfill
\begin{minipage}{0.46\textwidth}
\renewcommand{\arraystretch}{1.055}
\begin{tabular}{|l|l|}\hline
Identifier & Classification \\\hline\hline 
$\left[ 243, 16\right]$ & $ L(2)$ \\
$\left[ 243, 19\right]$ & $ Z'' \left( 3, 3 \right)$ \\
$\left[ 243, 20\right]$ & $ Z'\left( 3, 3 \right) $ \\
$\left[ 243, 24\right]$ & $ W \left( 1, 4 \right) $ \\
$\left[ 243, 25\right]$ & $ D\left( 0 \right) $ \\
$\left[ 243, 27\right]$ & $ X(9)$, $J(1)$ \\
$\left[ 243, 50\right]$ & $ Z \left( 3, 3 \right)$ \\
$\left[ 243, 55\right]$ & $U \left( 3, 2, 2 \right)$ \\
$\left[ 252, 11\right]$ & $ L_{7}^{(2)}(2,2) $ \\
$\left[ 279, 1 \right]$ & $ T_{31}^{(5)}(2)$ \\
$\left[ 300, 13\right]$ & $ \Delta \left( 6\times5^2,2 \right)$ \\
$\left[ 324, 3 \right]$ & $ \Delta \left( 3\times2^2,4 \right)$ \\
$\left[ 324, 13\right]$ & $ G \left( 1,2 \right) $ \\
$\left[ 324, 15\right]$ & $ Z \left( 3, 2, 2 \right)$ \\
$\left[ 324, 17\right]$ & $ Z' \left( 3, 2, 2 \right) $ \\
$\left[ 324, 43\right]$ & $ W \left( 2, 3 \right) $ \\
$\left[ 324, 45\right]$ & $Y \left( 2, 1 \right)$ \\
$\left[ 324, 49\right]$ & $ Z' \left( 6, 2 \right) $ \\
$\left[ 324, 51\right]$ & $ V(1)$ \\
$\left[ 324, 60\right]$ & $ Y(1)$ \\
$\left[ 324, 102\right]$ & $ \Delta ' \left( 6\times3^2,2,2 \right)$ \\
$\left[ 324, 111\right]$ & $ \Xi \left( 2,2 \right) $ \\
$\left[ 324, 128\right]$ & $ Z \left( 6, 2 \right)$ \\
$\left[ 333, 1 \right]$ & $ T_{37}^{(10)}(2) $ \\
$\left[ 351, 1 \right]$ & $ T_{13}^{(3)}(3)$ \\
$\left[ 351, 4 \right]$ & $ Q_{13}^{(3)}(2)$ \\
$\left[ 351, 5 \right]$ & $ Q_{13}^{(3)\prime}(2) $ \\
$\left[ 351, 7 \right]$ & $ P_{13}^{(3)}(2)$ \\
$\left[ 384, 571\right]$ & $ \Delta \left( 6\times4^2,3 \right)$ \\
$\left[ 384, 581\right]$ & $ S_4(5) $ \\
\hline
\end{tabular}
\caption{The finite subgroups of $U(3)$.
Part~2: groups with $243 \le \mbox{order} \le 384$.}
\label{tab:U3_par2} 
\end{minipage}
\end{centering}
\end{table}

\begin{table}
\begin{centering}
\begin{minipage}{0.46\textwidth}
\renewcommand{\arraystretch}{1.1}
\begin{tabular}{|l|l|}\hline
Identifier & Classification \\\hline\hline
$\left[ 387, 1 \right]$ & $ T_{43}^{(6)}(2)$ \\
$\left[ 432, 3 \right]$ & $ \Delta \left( 3\times4^2,3 \right)$ \\
$\left[ 432, 33\right]$ & $ \Delta \left( 6\times3^2,4 \right)$ \\
$\left[ 432, 57\right]$ & $ \Xi \left( 1,4 \right) $ \\
$\left[ 432, 100\right]$ & $ W \left( 4, 2 \right) $ \\
$\left[ 432, 102\right]$ & $ X(12)$ \\
$\left[ 432, 239\right]$ & $ \Pi \left( 1,2 \right) $ \\
$\left[ 432, 260\right]$ & $ \Delta \left( 6\times6^2,2 \right)$ \\
$\left[ 432, 273\right]$ & $ \hat \Xi \left( 1,3 \right) $ \\
$\left[ 441, 1 \right]$ & $ T_{49}^{(18)}(2) $ \\
$\left[ 441, 7 \right]$ & $ \Delta \left( 3\times7^2,2 \right)$ \\
$\left[ 468, 14\right]$ & $ L_{13}^{(3)} \left( 2,2 \right)$ \\
$\left[ 486, 26\right]$ & $ Z' \left( 3, 3, 1 \right) $ \\
$\left[ 486, 28\right]$ & $ Z \left( 3, 3, 1 \right)$ \\
$\left[ 486, 125\right]$ & $ H \left( 3,2,1 \right) $ \\
$\left[ 486, 164\right]$ & $ \Delta ' \left( 6\times3^2,3,1 \right)$ \\
$\left[ 513, 1 \right]$ & $ T_{19}^{(7)}(3)$ \\
$\left[ 513, 5 \right]$ & $ Q_{19}^{(7)}(2)$ \\
$\left[ 513, 6 \right]$ & $ Q_{19}^{(7)\prime}(2) $ \\
$\left[ 513, 8 \right]$ & $ P_{19}^{(7)}(2)$ \\
$\left[ 549, 1 \right]$ & $ T_{61}^{(13)}(2) $ \\ 
$\left[ 567, 1 \right]$ & $ T_{7}^{(2)}(4)$ \\
$\left[ 567, 4 \right]$ & $ Q_{7}^{(2)}(3)$ \\
$\left[ 567, 5 \right]$ & $ Q_{7}^{(2)\prime}(3)$ \\
$\left[ 567, 7 \right]$ & $ P_{7}^{(2)}(3)$ \\
$\left[ 567, 12\right]$ & $ S_{7}^{(2)\prime}(2)$ \\
$\left[ 567, 14\right]$ & $V_{7}^{(2)}(2)$ \\
$\left[ 567, 23\right]$ & $ Y_{7}^{(2)}(2) $ \\
$\left[ 567, 36\right]$ & $ S_{7}^{(2)}(2)$ \\
$\left[ 576, 3 \right]$ & $ \Delta \left( 3\times8^2,2 \right)$ \\
\hline
\end{tabular}
\caption{The finite subgroups of $U(3)$.
Part~3: groups with $387 \le \mbox{order} \le 576$.}
\label{tab:U3_par3} 
\end{minipage}
\hfill
\begin{minipage}{0.46\textwidth}
\renewcommand{\arraystretch}{1.13}
\begin{tabular}{|l|l|}\hline
Identifier & Classification \\\hline\hline
$\left[ 588, 16\right]$ & $ \Delta \left( 6\times7^2,2 \right)$ \\
$\left[ 600, 45\right]$ & $ \Delta \left( 6\times5^2,3 \right)$ \\
$\left[ 603, 1 \right]$ & $ T_{67}^{(29)}(2) $ \\
$\left[ 648, 19\right]$ & $ G \left( 1,3 \right) $ \\
$\left[ 648, 21\right]$ & $ Z \left( 3, 2, 3 \right)$ \\
$\left[ 648, 23\right]$ & $ Z' \left( 3, 2, 3 \right) $ \\
$\left[ 648, 244\right]$ & $ \Delta ' \left( 6\times3^2,2,3 \right)$ \\
$\left[ 648, 260\right]$ & $ Z \left( 6, 2, 1 \right)$ \\
$\left[ 648, 266\right]$ & $ \tilde Y(1) $ \\
$\left[ 648, 352\right]$ & $ \Xi \left( 2,3 \right) $ \\
$\left[ 648, 531\right]$ & $ \Upsilon(2) $ \\
$\left[ 648, 533\right]$ & $ \Omega(1)$ \\
$\left[ 648, 551\right]$ & $ \Theta(2)$ \\
$\left[ 648, 563\right]$ & $ \Delta ' \left( 6\times6^2,2,1 \right)$ \\
$\left[ 657, 1 \right]$ & $ T_{73}^{(8)}(2)$ \\
$\left[ 675, 5 \right]$ & $ \Delta \left( 3\times5^2,3 \right)$ \\
$\left[ 675, 9 \right]$ & $ W \left( 5, 2 \right) $ \\
$\left[ 675, 11\right]$ & $ X(15)$ \\
$\left[ 684, 11\right]$ & $ L_{19}^{(7)}(2,2)$ \\
$\left[ 711, 1 \right]$ & $ T_{79}^{(23)}(2) $ \\
$\left[ 729, 62\right]$ & $ L(3)$ \\
$\left[ 729, 63\right]$ & $ Z'' \left( 3, 4 \right)$ \\
$\left[ 729, 64\right]$ & $ Z' \left( 3, 4 \right) $ \\
$\left[ 729, 80\right]$ & $J(2) $               \\
$\left[ 729, 86\right]$ & $U \left( 3, 3, 2 \right)$ \\
$\left[ 729, 94\right]$ & $ W \left( 1, 5 \right) $ \\
$\left[ 729, 96\right]$ & $ \text{see section~\ref{U3sec5}} $ \\
$\left[ 729, 97\right]$ & $ \text{see section~\ref{U3sec5}} $ \\
$\left[ 729, 98\right]$ & $ \text{see section~\ref{U3sec5}} $ \\
$\left[ 729, 284\right]$ & $U \left( 3, 3, 3 \right)$ \\
\hline
\end{tabular}
\caption{The finite subgroups of $U(3)$.
Part~4: groups with $588 \le \mbox{order} \le 729$.}
\label{tab:U3_par4} 
\end{minipage}
\end{centering}
\end{table}

\begin{table}
\begin{centering}
\begin{minipage}{0.46\textwidth}
\renewcommand{\arraystretch}{1.1}
\begin{tabular}{|l|l|}\hline
Identifier & Classification \\\hline\hline
$\left[ 729, 393\right]$ & $ Z \left( 3, 4 \right)$ \\ 
$\left[ 729, 397\right]$ & $ Z \left( 9, 2 \right)$ \\
$\left[ 756, 11\right]$ & $ L_{7}^{(2)} \left( 2,3 \right) $ \\
$\left[ 756, 113\right]$ & $ M_{7}^{(2)}$ \\
$\left[ 756, 114\right]$ & $ M_{7}^{(2)\prime}$ \\
$\left[ 756, 116\right]$ & $ J_{7}^{(2)} $ \\
$\left[ 768, 1085333\right]$ & $ \Delta \left( 6\times4^2,4 \right)$ \\
$\left[ 768, 1085335\right]$ & $ \Delta \left( 6\times8^2,2 \right)$ \\
$\left[ 768, 1085351\right]$ & $ S_4(6) $ \\
$\left[ 819, 3 \right]$ & $ T_{91}^{(16)}(2) $ \\
$\left[ 819, 4 \right]$ & $ T_{91}^{(9)}(2)$ \\
$\left[ 837, 1 \right]$ & $ T_{31}^{(5)}(3)$ \\
$\left[ 837, 4 \right]$ & $ Q_{31}^{(5)\prime}(2) $ \\
$\left[ 837, 5 \right]$ & $ Q_{31}^{(5)}(2)$ \\
$\left[ 837, 7 \right]$ & $ P_{31}^{(5)}(2)$ \\
$\left[ 864, 69\right]$ & $ \Delta \left( 6\times3^2,5 \right)$ \\
$\left[ 864, 194\right]$ & $ \Xi \left( 1,5 \right) $ \\
$\left[ 864, 675\right]$ & $ \Pi \left( 1,3 \right) $ \\
$\left[ 864, 703\right]$ & $ \Delta \left( 6\times6^2,3 \right)$ \\
$\left[ 864, 737\right]$ & $ \hat \Xi \left( 1,4 \right) $ \\
$\left[ 873, 1 \right]$ & $ T_{97}^{(35)}(2) $ \\
$\left[ 900, 66\right]$ & $ \Delta \left( 3\times10^2,2 \right) $ \\
$\left[ 927, 1 \right]$ & $ T_{103}^{(46)}(2)$ \\
$\left[ 972, 3 \right]$ & $ \Delta \left( 3\times2^2,5 \right)$ \\
$\left[ 972, 29\right]$ & $ Z' \left( 3, 3, 2 \right) $ \\
$\left[ 972, 31\right]$ & $ Z \left( 3, 3, 2 \right)$ \\
$\left[ 972, 64\right]$ & $ \Delta \left( 6\times9^2,2 \right)$ \\
$\left[ 972, 117\right]$ & $ W \left( 2, 4 \right) $ \\
$\left[ 972, 121\right]$ & $ D(1) $ \\
\hline
\end{tabular}
\caption{The finite subgroups of $U(3)$.
Part~5: groups with $729 \le \mbox{order} \le 972$.}
\label{tab:U3_par5} 
\end{minipage}
\hfill
\begin{minipage}{0.46\textwidth}
\renewcommand{\arraystretch}{1.205}
\begin{tabular}{|l|l|}\hline
Identifier & Classification \\\hline\hline
$\left[ 972, 123\right]$ & $ X(18)$ \\
$\left[ 972, 147\right]$ & $Y \left( 3, 1 \right)$ \\
$\left[ 972, 152\right]$ & $ Z' \left( 6, 3 \right) $ \\
$\left[ 972, 153\right]$ & $ Z'' \left( 6, 3 \right)$ \\
$\left[ 972, 170\right]$ & see subsection~\ref{U3sec5} \\
$\left[ 972, 309\right]$ & $ H \left( 3,2,2 \right),\; G \left( 2,2 \right) $ \\
$\left[ 972, 348\right]$ & $ \Delta ' \left( 6\times3^2,3,2 \right)$ \\
$\left[ 972, 411\right]$ & $ \Xi \left( 3,2 \right) $ \\
$\left[ 972, 520\right]$ & $ Z \left( 6, 3 \right)$ \\
$\left[ 972, 550\right]$ & $U \left( 6, 2, 2 \right)$ \\
$\left[ 981, 1 \right]$ & $ T_{109}^{(45)}(2)$ \\
$\left[ 999, 1 \right]$ & $ T_{37}^{(10)}(3) $ \\
$\left[ 999, 5 \right]$ & $ Q_{37}^{(10)\prime}(2)$ \\
$\left[ 999, 6 \right]$ & $ Q_{37}^{(10)}(2) $ \\
$\left[ 999, 8 \right]$ & $ P_{37}^{(10)}(2) $ \\
$\left[ 1008, 57\right]$ & $ L_{7}^{(2)} \left( 4,2 \right) $ \\
$\left[ 1053, 16\right]$ & $ T_{13}^{(3)}(4)$ \\
$\left[ 1053, 25\right]$ & $ Q_{13}^{(3)\prime}(3) $ \\
$\left[ 1053, 26\right]$ & $ Q_{13}^{(3)}(3)$ \\
$\left[ 1053, 27\right]$ & $ P_{13}^{(3)}(3)$ \\
$\left[ 1053, 29\right]$ & $ Y_{13}^{(3)}(2)$ \\
$\left[ 1053, 32\right]$ & $ S_{13}^{(3)\prime}(2)$ \\
$\left[ 1053, 37\right]$ & $ V_{13}^{(3)}(2)$ \\
$\left[ 1053, 47\right]$ & $ S_{13}^{(3)}(2) $ \\
$\left[ 1089, 3\right]$ & $ \Delta \left( 3\times11^2,2 \right) $ \\
$\left[ 1116, 11\right]$ & $ L_{31}^{(5)} \left( 2,2 \right)$ \\
$\left[ 1143, 1\right]$ & $ T_{127}^{(19)}(2)$ \\
\hline
\end{tabular}
\caption{The finite subgroups of $U(3)$.
Part~6: groups with $972 \le \mbox{order} \le 1143$.}
\label{tab:U3_par6} 
\end{minipage}
\end{centering}
\end{table}

\begin{table}
\begin{centering}
\begin{minipage}{0.46\textwidth}
\renewcommand{\arraystretch}{1.203}
\begin{tabular}{|l|l|}\hline
Identifier & Classification \\\hline\hline
$\left[ 1161, 6\right]$ & $ T_{43}^{(6)}(3)$ \\
$\left[ 1161, 10\right]$ & $ Q_{43}^{(6)}(2)$ \\
$\left[ 1161, 11\right]$ & $ Q_{43}^{(6)\prime}(2) $ \\
$\left[ 1161, 12\right]$ & $ P_{43}^{(6)}(2)$ \\
$\left[ 1176, 57\right]$ & $ \Delta \left( 6\times7^2,3 \right)$ \\
$\left[ 1197, 3\right]$ & $ T_{133}^{(11)}(2)$ \\
$\left[ 1197, 4\right]$ & $ T_{133}^{(30)}(2)$ \\
$\left[ 1200, 183 \right]$ & $ \Delta \left( 6\times5^2,4 \right)$ \\
$\left[ 1200, 682 \right]$ & $ \Delta \left( 6\times10^2,2 \right) $ \\
$\left[ 1251, 1\right]$ & $ T_{139}^{(42)}(2)$ \\
$\left[ 1296, 3\right]$ & $ \Delta \left( 3\times4^2,4 \right)$ \\
$\left[ 1296, 35\right]$ & $ G \left( 1,4 \right) $ \\
$\left[ 1296, 37\right]$ & $ Z \left( 3, 2, 4 \right)$ \\
$\left[ 1296, 39\right]$ & $ Z' \left( 3, 2, 4 \right) $ \\
$\left[ 1296, 220 \right]$ & $ W \left( 4, 3 \right) $ \\
$\left[ 1296, 222 \right]$ & $Y \left( 2, 2 \right)$ \\
$\left[ 1296, 226 \right]$ & $ V(2)$ \\
$\left[ 1296, 227 \right]$ & $ Z' \left( 12, 2 \right)$ \\
$\left[ 1296, 237 \right]$ & $ Y(2)$ \\
$\left[ 1296, 647 \right]$ & $ \Delta ' \left( 6\times3^2,2,4 \right)$ \\
$\left[ 1296, 688 \right]$ & $ Z' \left( 6, 2, 2 \right) $ \\
$\left[ 1296, 689 \right]$ & $ Z \left( 6, 2, 2 \right)$ \\
$\left[ 1296, 699 \right]$ & see subsection~\ref{U3sec4} \\
$\left[ 1296, 1239\right]$ & $ \Xi \left( 2,4 \right) $ \\
$\left[ 1296, 1499\right]$ & $ Z \left( 12, 2 \right) $ \\
$\left[ 1296, 1995\right]$ & $ \Pi \left( 2,2 \right) $ \\
$\left[ 1296, 2113\right]$ & $ \Delta ' \left( 6\times6^2,2,2 \right)$ \\
$\left[ 1296, 2203\right]$ & $ \hat \Xi \left( 2,3 \right) $ \\
\hline
\end{tabular}
\caption{The finite subgroups of $U(3)$.
Part~7: groups with $1161 \le \mbox{order} \le 1296$.}
\label{tab:U3_par7} 
\end{minipage}
\hfill
\begin{minipage}{0.46\textwidth}
\renewcommand{\arraystretch}{1.2}
\begin{tabular}{|l|l|}\hline
Identifier & Classification \\\hline\hline
$\left[ 1323, 1\right]$ & $ T_{49}^{(18)}(3) $ \\
$\left[ 1323, 4\right]$ & $ Q_{49}^{(18)}(2) $ \\
$\left[ 1323, 5\right]$ & $ Q_{49}^{(18)\prime}(2)$ \\
$\left[ 1323, 7\right]$ & $ P_{49}^{(18)}(2) $ \\
$\left[ 1323, 14\right]$ & $ \Delta \left( 3\times7^2,3 \right)$ \\
$\left[ 1323, 40\right]$ & $W \left( 7, 2 \right) $ \\
$\left[ 1323, 42\right]$ & $ X(21)$ \\
$\left[ 1332, 14\right]$ & $ L_{37}^{(10)} \left( 2,2 \right) $ \\
$\left[ 1359, 1\right]$ & $ T_{151}^{(32)}(2)$ \\
$\left[ 1404, 14\right]$ & $ L_{13}^{(3)} \left( 2,3 \right)$ \\
$\left[ 1404, 137 \right]$ & $ M_{13}^{(3)}$ \\
$\left[ 1404, 138 \right]$ & $ M_{13}^{(3)\prime}$ \\
$\left[ 1404, 140 \right]$ & $ J_{13}^{(3)}$ \\
$\left[ 1413, 1\right]$ & $ T_{157}^{(12)}(2)$ \\
$\left[ 1452, 11\right]$ & $ \Delta \left( 6\times11^2,2 \right) $ \\
$\left[ 1458, 615 \right]$ & $ Z' \left( 3, 4, 1 \right) $ \\
$\left[ 1458, 618 \right]$ & $ Z \left( 3, 4, 1 \right)$ \\
$\left[ 1458, 663 \right]$ & $ \text{see section~\ref{U3sec5}} $ \\
$\left[ 1458, 666 \right]$ & $ \text{see section~\ref{U3sec5}} $ \\
$\left[ 1458, 1095\right]$ & $ H \left( 3,3,1 \right) $ \\
$\left[ 1458, 1354\right]$ & $ \Delta ' \left( 6\times3^2,4,1 \right)$ \\
$\left[ 1458, 1371\right]$ & $ \Delta ' \left( 6\times9^2,2,1 \right)$ \\
$\left[ 1467, 1\right]$ & $ T_{163}^{(58)}(2)$ \\
$\left[ 1521, 1\right]$ & $ T_{169}^{(22)}(2)$ \\
$\left[ 1521, 7\right]$ & $ \Delta \left( 3\times13^2,2 \right) $ \\
$\left[ 1536, 408544641\right]$ & $ \Delta \left( 6\times8^2,3 \right)$ \\
$\left[ 1536, 408544678\right]$ & $ \Delta \left( 6\times4^2,5 \right)$ \\
$\left[ 1536, 408544687\right]$ & $ S_4(7) $ \\
\hline
\end{tabular}
\caption{The finite subgroups of $U(3)$.
Part~8: groups with $1323 \le \mbox{order} \le 1536$.}
\label{tab:U3_par8} 
\end{minipage}
\end{centering}
\end{table}

\begin{table}
\begin{centering}
\begin{minipage}{0.46\textwidth}
\renewcommand{\arraystretch}{1.2}
\begin{tabular}{|l|l|}\hline
Identifier & Classification \\\hline\hline
$\left[ 1539, 16\right]$ & $ T_{19}^{(7)}(4)$ \\
$\left[ 1539, 25\right]$ & $ Q_{19}^{(7)\prime}(3) $ \\
$\left[ 1539, 26\right]$ & $ Q_{19}^{(7)}(3)$ \\
$\left[ 1539, 27\right]$ & $ P_{19}^{(7)}(3)$ \\
$\left[ 1539, 29\right]$ & $ Y_{19}^{(7)}(2)$ \\
$\left[ 1539, 32\right]$ & $ S_{19}^{(7)\prime}(2)$ \\
$\left[ 1539, 37\right]$ & $V_{19}^{(7)}(2)$ \\
$\left[ 1539, 47\right]$ & $ S_{19}^{(7)}(2) $ \\
$\left[ 1548, 11\right]$ & $ L_{43}^{(6)} \left( 2,2 \right)$ \\
$\left[ 1575, 7\right]$ & $ L_{7}^{(2)} \left( 5,2 \right) $ \\
$\left[ 1629, 1\right]$ & $ T_{181}^{(48)}(2)$ \\
$\left[ 1647, 6\right]$ & $ T_{61}^{(13)}(3) $ \\
$\left[ 1647, 10\right]$ & $ Q_{61}^{(13)}(2) $ \\
$\left[ 1647, 11\right]$ & $ Q_{61}^{(13)\prime}(2)$ \\
$\left[ 1647, 12\right]$ & $ P_{61}^{(13)}(2) $ \\
$\left[ 1701, 68\right]$ & $ T_{7}^{(2)}(5)$ \\
$\left[ 1701, 102 \right]$ & $ \text{see section~\ref{U3sec4}} $ \\
$\left[ 1701, 112 \right]$ & $ \text{see section~\ref{U3sec5}} $ \\
$\left[ 1701, 115 \right]$ & $ S_{7}^{(2)\prime}(3)$ \\
$\left[ 1701, 126 \right]$ & $ Q_{7}^{(2)\prime}(4)$ \\
$\left[ 1701, 127 \right]$ & $ Q_{7}^{(2)}(4)$ \\
$\left[ 1701, 128 \right]$ & $ P_{7}^{(2)}(4)$ \\
$\left[ 1701, 130 \right]$ & $ \text{see section~\ref{U3sec5}} $ \\
$\left[ 1701, 131 \right]$ & $ \text{see section~\ref{U3sec5}} $ \\
$\left[ 1701, 138 \right]$ &$V_7^{(2)}(3)$ \\
$\left[ 1701, 240 \right]$ & $ S_{7}^{(2)}(3)$ \\
$\left[ 1701, 261 \right]$ & $ Y_{7}^{(2)}(3)$ \\
\hline
\end{tabular}
\caption{The finite subgroups of $U(3)$.
Part~9: groups with $1539 \le \mbox{order} \le 1701$.}
\label{tab:U3_par9} 
\end{minipage}
\hfill
\begin{minipage}{0.46\textwidth}
\renewcommand{\arraystretch}{1.2}
\begin{tabular}{|l|l|}\hline
Identifier & Classification \\\hline\hline
$\left[ 1728, 3\right]$ & $ \Delta \left( 3\times8^2,3 \right)$ \\
$\left[ 1728, 185 \right]$ & $ \Delta \left( 6\times3^2,6 \right)$ \\
$\left[ 1728, 953 \right]$ & $ \Xi \left( 1,6 \right) $ \\
$\left[ 1728, 1286\right]$ & $W \left( 8, 2 \right) $ \\
$\left[ 1728, 1290\right]$ & $ X(24)$ \\
$\left[ 1728, 2785\right]$ & $ \Pi \left( 1,4 \right) $ \\
$\left[ 1728, 2847\right]$ & $ \Delta \left( 6\times12^2,2 \right) $ \\
$\left[ 1728, 2855\right]$ & $ \Delta \left( 6\times6^2,4 \right)$ \\
$\left[ 1728, 2929\right]$ & $ \hat \Xi \left( 1,5 \right) $ \\
$\left[ 1737, 1\right]$ & $ T_{193}^{(84)}(2)$ \\
$\left[ 1764, 11\right]$ & $ L_{49}^{(18)} \left( 2,2 \right) $ \\
$\left[ 1764, 91\right]$ & $ \Delta \left( 3\times14^2,2 \right) $ \\
$\left[ 1791, 1\right]$ & $ T_{199}^{(92)}(2)$ \\
$\left[ 1809, 6\right]$ & $ T_{67}^{(29)}(3) $ \\
$\left[ 1809, 10\right]$ & $ Q_{67}^{(29)}(2) $ \\
$\left[ 1809, 11\right]$ & $ Q_{67}^{(29)\prime}(2)$ \\
$\left[ 1809, 12\right]$ & $ P_{67}^{(29)}(2) $ \\
$\left[ 1872, 60\right]$ & $ L_{13}^{(3)} \left( 4,2 \right)$ \\
$\left[ 1899, 1\right]$ & $ T_{211}^{(14)}(2)$ \\
$\left[ 1944, 35\right]$ & $ Z' \left( 3, 3, 3 \right) $ \\
$\left[ 1944, 37\right]$ & $ Z \left( 3, 3, 3 \right)$ \\
$\left[ 1944, 70\right]$ & $ \Delta \left( 6\times9^2,3 \right)$ \\
$\left[ 1944, 707 \right]$ &
$ H \left( 3,2,3 \right),\; G \left( 2,3 \right) $ \\
$\left[ 1944, 746 \right]$ & $ \Delta ' \left( 6\times3^2,3,3 \right)$ \\
$\left[ 1944, 832 \right]$ & $ Z' \left( 6, 3, 1 \right) $ \\
$\left[ 1944, 833 \right]$ & $ Z \left( 6, 3, 1 \right)$ \\
$\left[ 1944, 1123\right]$ & $ \Xi \left( 3,3 \right) $ \\
\hline
\end{tabular}
\caption{The finite subgroups of $U(3)$.
Part~10: groups with $1728 \le \mbox{order} \le 1944$.}
\label{tab:U3_par10} 
\end{minipage}
\end{centering}
\end{table}

\begin{table}
\begin{centering}
\begin{minipage}{0.46\textwidth}
\renewcommand{\arraystretch}{1.2}
\begin{tabular}{|l|l|}\hline
Identifier & Classification \\\hline\hline
$\left[ 1944, 2293\right]$ & $ \Upsilon(3) $ \\
$\left[ 1944, 2294\right]$ & $ \Upsilon'(3)$ \\
$\left[ 1944, 2333\right]$ & $ \Theta(3)$ \\
$\left[ 1944, 2363\right]$ & $ H \left( 6,2,1 \right) $ \\
$\left[ 1944, 2415\right]$ & $ \Delta ' \left( 6\times6^2,3,1 \right)$ \\
$\left[ 1944, 3448\right]$ & $ \Omega(2)$ \\
$\left[ 1953, 3\right]$ & $ T_{217}^{(25)}(2)$ \\
$\left[ 1953, 4\right]$ & $ T_{217}^{(67)}(2)$ \\
$\left[ 1971, 6\right]$ & $ T_{73}^{(8)}(3)$ \\
$\left[ 1971, 10\right]$ & $ Q_{73}^{(8)\prime}(2) $ \\
$\left[ 1971, 11\right]$ & $ Q_{73}^{(8)}(2)$ \\
$\left[ 1971, 12\right]$ & $ P_{73}^{(8)}(2)$ \\
\hline
\end{tabular}
\caption{The finite subgroups of $U(3)$.
Part~11: groups with $1944 \le \mbox{order} < 2000$.}
\label{tab:U3_par11} 
\end{minipage}
\end{centering}
\end{table}

\end{document}